\documentclass[german
              ,english
              ,final
              ,twoside
              ,11pt
              ]{amsart}
\usepackage{gt-sim}
\usepackage[latin1]{inputenc}
\usepackage{babel}
\usepackage{cas-math}
\usepackage{cas-paper}
\usepackage[arrow,cmtip,matrix]{xy}

\newcommand\draftcomment[1]{{}}

\newcommand\A{{\mathcal A}}
\newcommand\mintimes{\mathop{\hat\times}}
\newcommand\dpairp[3]{(\Delta #1_{[#2,#3]},\Delta #1_{[#2,#3)})}
\newcommand\pairp[3]{(#1_{[#2,#3]},#1_{[#2,#3)})}
\newcommand\dpair[3]{\Delta #1_{[#2,#3]},\Delta #1_{[#2,#3)}}
\newcommand\Hp{{}'\!H}
\newcommand\Hpp{{}''\!H}

\foreignlanguage{german}{
\hyphenation{ar-range-ment ar-range-ments MacPher-son Arnim-al-lee}}

\begin{document}

\title[The cohomology ring of a projective arrangement]
{The cohomology ring of the complement of a finite family of linear subspaces
  in a complex projective space}

\author{Carsten Schultz} 
\address{\foreignlanguage{german}{Carsten Schultz\\Institut f\"ur
Mathematik II, Freie Universit\"at Berlin, Arnimallee 2--6, D-14195 Berlin}}
\email{carsten@codimi.de}
\date{11th November 2004}

\begin{abstract}
The integral cohomology ring of the complement of an arrangement of
linear subspaces of a finite dimensional complex projective space is
determined by combinatorial data, i.e.\ the intersection poset and the
dimension function.
\end{abstract}

\subjclass[2000]{52C35, 55N45; 05E25}
\keywords{subspace arrangement,
          projective arrangement,
          cohomology ring, 
          Orlik-Solomon algebra
         }

\maketitle

\section{Introduction}
\label{sec:intro}
Let $V$ be an $(n+1)$--dimensional complex vector space and $\A$ a
linear arrangement in $V$, i.e.\ a finite set of proper linear
subspaces of $V$.  We define $Q$ to be the partially ordered set
(poset) of intersections of $\A$,
that is the set $\set{\Intersection C\colon C\subset\A}$ ordered by
inclusion.  $Q$~is called the \emph{intersection poset} of~$\A$.
Note that $Q$ has the maximum $V=\Intersection\emptyset$.

A recurrent theme in the study of arrangements are questions of the
type if a certain property of an arrangement is determined by
combinatorial data, and if so, how it can be described using only
those.  For example, De~Concini and Procesi have shown that the
intersection poset labelled with the dimensions of the intersections
determines the rational cohomology ring
$H^\ast\left(V\wo\Union\A;\Q\right)$ of the complement of the
arrangement~\cite{conciniprocesi95}.  Yuzvinsky derived an explicit
formula from their work~\cite{yuzvinsky97}.  Deligne, Goresky, and
MacPherson and indepently Mark de Longueville and the current author
extended this result to cohomology with integral coefficients and
clarified the extent to which it holds for real
arrangements~\cite{delgormacph,delong-schultz}.

We will generalize this result to the complement of the projective
arrangement $P\A\deq\set{PA\colon A\in\A}$ in the $n$-dimensional
complex projective space $PV$.  This will be a generalization, because
the complement of an affine arrangement is the same as the complement
of a projective arrangement containing an additional hyperplane,
details will be given in \prettyref{rem:affine-homology} and
\prettyref{rem:affine-products}.  It will not imply the results on
real real linear arrangements mentioned above, though.

\begin{nota}
For $u\in Q$, we set
\[d(u)\deq\dim_\C(u)-1=\dim_\C(Pu).\]
We will find it useful to write
\[Q_{[r,s)}\deq\set{u\in Q\colon r\le d(u)\lt s}\]
etc.~for certain subsets of $Q$.
\end{nota}

Our goal is thus to describe the cohomology ring $H^\ast\left(PV\wo
\Union P\A\right)$ in terms of the intersection poset~$Q$ and the dimension
function~$d$.

\begin{nota}
For a poset $P$, $\Delta P$ denotes the order complex of $P$, i.e.\ the
simplicial complex with $P$ as the set of vertices and the chains in
$P$ as simplexes.
\end{nota}

By Poincar\'e duality, 
$H^i(PV\wo \Union P\A)\isom H_{2n-i}(PV, \Union P\A)$, and the
homology of $(PV, \Union P\A)$ has been given by Goresky an MacPherson 
\cite[Thm.~D]{goreskymacpherson} as
\begin{equation}\label{eqn:intro-gm}
H_i(PV, {\textstyle\Union P\A}) \isom
    \dsum_{k=0}^{\lfloor i/2\rfloor} 
     H_{i-2k}
      \left(\dpair Qkn\right).
\end{equation}
Actually, they state the cohomology analogue of this formula.

The proof we give in \prettyref{prop:proj-homology} constructs the
isomorphism using explicit homomorphisms
\begin{equation}\label{eqn:intro-hk}
h_k\colon H_{\bullet-2k}
      \left(\dpair Qkn\right)
 \to H_\bullet(PV, {\textstyle\Union P\A})
\end{equation}
(see \prettyref{def:hk})
in the spirit of the analogous construction for linear arrangements
employed by Ziegler and \v{Z}ivaljevi\'c~\cite{zieglerzivaljevic}.

Our goal now is to describe cup products in $H^\ast(PV\wo\Union P\A)$,
respectively the corresponding intersection products in
$H_\ast(PV,\Union P\A)$.  For the formulation of the corresponding
result for linear arrangements in \cite{delong-schultz}, the following
product on the homology of an order complex is central.

\begin{defn}
Noting that the map
\begin{align*}
\lmin\colon Q\times Q&\to Q\\
u\lmin v&\deq u\intersect v
\intertext{is a simplicial map $\Delta Q\times\Delta Q\to\Delta
Q$, we define}
\mintimes\colon C_\ast(\Delta Q)\tens C_\ast(\Delta Q)
    &\to C_\ast(\Delta Q)\\
 c\mintimes d&\deq\lmin_\ast(c\times d)    
\end{align*}
and denote the corresponding operation on homology by the same symbol.
\end{defn}

With this notation, our main result 
(\prettyref{thm:main}) can be formulated:  
For $c\in H_\ast
      \left(\dpair Qkn\right)$, $d\in H_\ast
      \left(\dpair Qln\right)$,
\begin{equation}\label{eqn:intro-product}
    h_k(c)\iprod h_l(d)
    =\begin{cases}
         h_{k+l-n}(c\mintimes d),& k+l\ge n.\\
	 0,& k+l<n.
     \end{cases}
\end{equation}
While our setting is very similar to that in \cite{delong-schultz},
the proof cannot be as direct as the one for linear arrangements
presented there, as additional problems arise.  In particular, in the
case $k+l<n$ of \prettyref{eqn:intro-product} we will not be able to
arrange things in such a way that the chains which will naturally
represent $h_k(c)$ and $h_l(d)$ do not intersect geometrically, or
only in $\Union P\A$.  This is in contrast to the corresponding case
for linear arrangements (Prop.~4.3 in \cite{delong-schultz}). 
One step towards the
solution will be to not only consider the order complex $\Delta Q$,
but also, as in \cite{delgormacph}, to consider the the poset $Q$ itself
to be a space, namely the space of strata, equipped with the quotient
topology.  The combination of both will prove useful in 
\prettyref{prop:detect-h0}.

Finally, in \prettyref{sec:c-arr} we specialize the result to the case
of projective $c$-arrangements obtaining a presentation of the
cohomology ring of the complement in terms of generators and relations
in the spirit of Orlik and Solomon.  This complements results on
linear arrangements by Feichtner and Ziegler~\cite{feichtnerziegler}
and proceeds in the same way as Yuzvinskys derivation of a
presentation of the cohomology ring of a linear arrangement with
geometric intersection poset~\cite{yuzvinsky99}.

This article is a presentation of the main result of the author's
thesis\footnote{November 2004:  The thesis has been submitted, but not
yet defended.}~\cite{coa}.  The thesis
describes the spectral sequence at the core
of~\prettyref{sec:homology} and \prettyref{sec:graded-ring} in a less
ad hoc manner, considers real arrangements in more detail, and
inspects the relationship between affine and projective
arrangements more closely.

I wish to thank my thesis advisor Elmar Vogt for the support during
the work on the results presented here.  They are a follow-up to joint
work with Mark de Longueville.  Working with him was a very pleasant
experience, and I want to thank him for also always being open to
discussions of the current work.

\section{The homology spectral sequence and the direct sum
  decomposition of the homology of the arrangement}
\label{sec:homology}
In this section we construct the isomorphism~\prettyref{eqn:intro-gm}.
\begin{defn} \label{def:X}
We define a bigraded abelian group~$X$ by
\begin{equation*}
X_{p,q}\deq\dsum_{\substack{u\in Q^{q+1}\\u_0<\dots<u_q=V}} S_p(Pu_0),
\end{equation*}
where $S_\ast$ denotes the simplicial chain complex.  
Writing
$c\tens\simp{u_0,\dots,u_q}$ for the element
$c\in S_\ast(Pu_0)$ of the summand indexed by
$(u_0,\dots,u_q)$, we make this into a double chain complex by
\begin{align*}
\dd'\colon X_{p,q}&\to X_{p-1,q}\\
c\tens\simp{u_0,\dots,u_q}&\mapsto\dd c\tens\simp{u_0,\dots,u_q}\\
\intertext{and}
	\dd''\colon X_{p,q}&\to X_{p,q-1}\\
	c\tens\simp{u_0,\dots,u_q}&
	\mapsto\sum_{i=0}^{q-1}(-1)^{p+i}c
		\tens\simp{u_0,\dots,\hat{u_i},\dots,u_q},
\end{align*}
satisfying $\dd'\dd''+\dd''\dd'=0$.  Note that in the definition of
$\dd''$ we regard $S_\ast(u_0)$ as a subcomplex of $S_\ast(u_1)$.
\end{defn}

The rough idea of the following is that on one hand
$\Hp(\Hpp_q(X))\isom0$ for $q\ne0$ and 
$\Hp(\Hpp_0(X))\isom H(PV,\Union P\A)$, while on the other hand
$\Hpp_{2k}(\Hp_p(X))\isom H_{2k}(\C P^k)\tens H_q(\dpair Qkn)$, 
$\Hpp_{2k+1}(\Hp_p(X))\isom 0$,
so that
there exists a spectral sequence with $E^2$--term isomorphic to the
latter that converges against $H(PV,\Union P\A)$.  Convergence alone
will not suffice, though.  Not only will the $E^2$--term equal the
$E^\infty$--term, all extensions will be trivial.

\begin{prop}\label{prop:X}
The map
\begin{align*}
\eps\colon X_{p,q}&\to S_p(PV,\Union P\A)\\
c\tens\simp{u_0,\dots,u_q}&\mapsto0,\qquad\text{if $q>0$,}\\
c\tens\simp{V}&\mapsto c
\end{align*}
is a chain map from the total complex of X to the relative simplicial
chain complex and induces an isomorphism
$\eps_\ast\colon H_\ast(X)\xto{\isom}H_\ast(PV,\Union P\A)$.
\end{prop}

\begin{proof}
For a simplicial simplex $\sigma\in S_\ast(PV)$ the subset
$\set{u\colon \im\sigma\subset Pu}$ of $Q$ has a minimum that we call 
$u_\sigma$.
We define
\begin{align*}
    K\colon X_{p,q}&\to X_{p,q+1} \\
    \sigma\tens\simp{u_0,\dots,u_q}
        &\mapsto\begin{cases}
                    (-1)^p\sigma\tens\simp{u_\sigma,u_0,\dots,u_q},
                      &u_0>u_\sigma,\\
		    0,&u_0=u_\sigma,
                \end{cases}
\end{align*}
and (do not have to) calculate $(\dd''K+K\dd'')x = x$
for $x\in X_{p,q}$, $q>0$,
so that $\Hpp_q(X_{p,\bullet})=0$ for $q>0$.
Furthermore
$\im \left(X_{p,1}\xto{\dd''}X_{p,0}\right)
 =\left(\sum_{u\in Q_{[0,n)}} S_p(Pu)\right)\tens\simp V$.
The claim now follows from the homology version of
\cite[Thm~I.4.8.1]{godement}
and the fact that the inclusion 
$\sum_{u\in Q_{[0,n)}} S_p(Pu)\to S_\ast(\Union P\A)$ induces an
isomorphism in homology.
\end{proof}
We start the construction of the homomorphisms~$h_k$ announced
in~\prettyref{eqn:intro-hk}.
\begin{prop}
Let $k\in\N$.  We consider all functions $x$ assigning to each $u\in
Q_{[k,n]}$ a system $(x^u_j)_{0\le j\le k}$ of $k+1$ vectors in $u$.
We regard the set of all these functions as the complex affine space
$\prod_{u\in Q_{[k,n]}}u^{k+1}$.  Among those, the set of all
functions such that for all $u_0,\dots,u_r\in Q_{[k,n]}$,
$\lambda_0,\dots,\lambda_r\in\C$ with $u_0<u_1<\dots<u_r$ and
$(\lambda_0,\dots,\lambda_r)\ne0$, the system of vectors
\begin{equation*}
\left(\sum_{i=0}^r \lambda_i x^{u_i}_j\right)_{0\le j\le k}
\end{equation*}
is linearly independent
form a non-empty, Zariski-open set.
\end{prop}

\begin{proof}
Since a finite intersection of non-empty Zariski-open sets is again
non-empty and Zariski-open, it suffices to consider a fixed chain 
$u_0<u_1<\dots<u_r$, $u_i\in Q_{[k,n]}$.  The set
\begin{multline*}
\Biggl\{
    (x^0,\dots,x^r)\in u_0^{k+1}\times\dots\times u_r^{k+1}
    \colon\\
    \text{There exist $\lambda\in\C^{r+1}$, $\mu\in\C^{k+1}$,
    $\lambda,\mu\ne0$
    with\ }
    \sum_{j=0}^k\sum_{i=0}^r\mu_j\lambda_i x^i_j=0
\Biggr\}
\end{multline*}
is algebraic by the 
main theorem of elimination theory\cite[I.5,~thm~3]{shafarevich1}, 
because the
equation is homogenous in $\lambda$ and $\mu$.
To see that the complement of this set is non-empty, we choose a basis 
$(e_l)_{l=0,\dots,n}$ of $V$ such that $e_l\in u_i$ for $l\le i+k$
and set $x^i_j\deq e_{i+j}$.  Now if 
$\sum_j\sum_i\mu_j\lambda_i x^i_j=0$ then
$\sum_i \lambda_i\mu_{s-i}=0$ for all~$s$, and it follows that
$\lambda=0$, $\mu=0$.
\end{proof}

\begin{defprop}  \label{def:fk}
Let $k\in\N$ and $x$ be any function as in the preceding proposition,
but linear independence of $\left(\sum_{i=0}^r \lambda_i
x^{u_i}_j\right)$ only required for $\lambda_i\ge0$,
$\sum\lambda_i=1$.  Writing points of $\simp{u_0,\dots,u_r}$ as
$\sum_{i=0}^r \lambda_i u_i$, we define a map
\begin{align*}
f^k\colon \C P^k\times\Delta Q_{[k,n]}&\to PV\\
\left([\mu_0:\dots:\mu_k],\sum_{i=0}^r \lambda_i u_i\right)
	&\mapsto
		\left[\sum_{j=0}^k\sum_{i=0}^r\mu_j\lambda_i x^{u_i}_j
		\right].
\end{align*}
The homotopy class of $f_k$ does not depend on the particular choice
of $x$.
\end{defprop}

\begin{proof}
The non-dependence of the homotopy class of~$f_k$ on the choice of~$x$
follows from the space of possibly choices, which contains a
Zariski-open set, being path connected.
\end{proof}

\begin{defn}
We define another, albeit degenerated, double complex by
\begin{equation*}
Y_{p,q}\deq\begin{cases}
	C_q(\Delta Q_{[k,n]}, \Delta Q_{[k,n)}),&p=2k,\\
	0,&p=2k+1,
\end{cases}
\end{equation*}
$\dd'=0$ and $\dd''$ being the simplicial boundary operator.
\end{defn}

\begin{prop}\label{prop:X-Y}
Let $o_k\in S_{2k}(\C P^k)$ denote a chain representing the
orientation class.  The map
\begin{align*}
g\colon Y_{2k,q}&\to X^{\rm tot}_{2k+q}\\
\simp{u_0,\dots,u_q}&\mapsto
	\sum_{i=0}^q f^k_\ast(o_k\times\simp{u_0,\dots,u_i})
		\tens\simp{u_i,\dots,u_q}
\end{align*}
is a chain map between the total complexes and induces an isomorphism
$H_\ast(Y)\to H_\ast(X)$.
\end{prop}

\begin{proof}
The fact that $g$ is a chain map follows from the usual calculation
used to verify this property for the Alexander-Whitney diagonal
approximation.

$g[Y_{\ast,q}]\subset\dsum_{q'\le q}X_{\ast,q'}$, that is $g$ respects
the filtrations by $q$ and therefore induces a homomorphism between
the spectral sequences defined by those filtrations.  The induced map
between the $E^1$--terms $\Hp_p(Y_{\bullet,q})=Y_{p,q}$ and
$\Hp_p(X_{\bullet,q})=\dsum_{u_0<\dots<u_q=V} H_p(Pu_0)$ 
is given, for $p=2k$, by
\begin{equation*}
   \simp{u_0,\dots,u_q}\mapsto
	[f^k_\ast(o_k\times\simp{u_0})]
		\tens\simp{u_0,\dots,u_q}.
\end{equation*}
$f^k(\cdot, u_0)\colon\C P^k\to Pu_0$ is a linear embedding and
therefore $[f^k_\ast(o_k\times\simp{u_0})]$ a generator of
$H_p(Pu_0)$.  So $g$ induces an isomorphism between the $E^1$--terms
and 
by \cite[Thm~I.4.3.1]{godement} it follows that 
the induced map between
the homology groups of the total complexes is also an isomorphism.
\end{proof}

\begin{defn}   \label{def:hk}
We define a homomorphism
\begin{align*}
h_k\colon H_r(\Delta Q_{[k,n]}, \Delta Q_{[k,n)})
	&\to H_{2k+r}(PV, \Union P\A) \\
c &\mapsto f^k_\ast([\C P^k]\times c)
\end{align*}
where $[\C P^k]$ denotes the orientation class of $\C P^k$.
\end{defn}

We can now prove the main result of this section.

\begin{prop}  \label{prop:proj-homology}
The map
\begin{equation*}
\sum_{k=0}^n h_k\colon\dsum_{k=0}^n 
             H_\ast(\Delta Q_{[k,n]}, \Delta Q_{[k,n)})[-2k]
	\to H_\ast(PV, \Union P\A)
\end{equation*}
is an isomorphism of graded abelian groups, in particular
\begin{align*}
H^{2n-r}(PV\wo\Union P\A)
 &\isom H_r(PV, \Union P\A)\\
 &\isom \dsum_k H_{r-2k}(\Delta Q_{[k,n]}, \Delta Q_{[k,n)})\\
 &\isom \dsum_k \tilde H_{r-2k-1}(\Delta Q_{[k,n)}).
\end{align*}
\end{prop}

\begin{proof}
The map $\sum h_k$ is just the composition of the isomorphisms
considered in \prettyref{prop:X-Y} and \prettyref{prop:X}.

For the isomorphism $ H_{s}(\Delta Q_{[k,n]}, \Delta
Q_{[k,n)})\xto[\isom]\dd \tilde H_{s-1}(\Delta Q_{[k,n)})$, note that
$\Delta Q_{[k,n]}$ is a cone and hence acyclic, because $ Q_{[k,n]}$
has a maximal element,~$V$.
\end{proof}

\begin{rem}  \label{rem:affine-homology}
If there is an $A_0\in\A$ with $d(A_0)=n-1$, i.e.\ if the arrangement
contains a hyperplane, then the complement of $\A$ can be considered
the complement of an affine arrangement in the affine space $PV\wo
PA_0$.  Since every affine arrangement arises in this way by adding a
hyperplane at infinity, the usual formula for the complement of an
affine (or linear) arrangement can be derived from the formula for
projective arrangements.  This is called an `interesting exercise' in
\cite{goreskymacpherson}.

We sketch how to do this.  Set $Q'\deq Q\wo\set{A_0\lmin q\colon q\in
Q}$.  $Q'$ is the poset of all intersections that are not contained in
$A_0$ and is can be defined as the intersection poset of the affine
arrangement in question, without reference to~$A_0$. For $0\le k<n$,
the simplical complex $\Delta(Q_{[k,n)}\wo Q'_{\set k})$ is acyclic,
since it contains the cone $\Delta\set{q\colon q\le A_0, d(q)\ge k}$
as a deformation retract.  Therefore the first map in
\begin{align*}
H\bigl(\dpair Qkn\bigr)
&\xto\isom 
  H\bigl(\Delta Q_{[k,n]},
    \Delta Q_{[k,n)}\unite\Delta(Q_{[k,n]}\wo Q'_{\set k})\bigr)
\\&\xot\isom
H\bigl(\Delta Q'_{[k,n]},
       \Delta Q'_{[k,n)}\unite \Delta Q'_{(k,n]}\bigr),
\end{align*}
which is induced by inclusion, is an isomorphism.  The second map is
also induced by inclusion and is an isomorphism by excision.
The last group splits as
\begin{equation*}
\dsum_{u\in Q'_{\set k}} H(\Delta[u,V],\Delta[u,V)\unite\Delta(u,V))
 \xto\isom
H\bigl(\Delta Q'_{[k,n]},
       \Delta Q'_{[k,n)}\unite \Delta Q'_{(k,n]}\bigr),
\end{equation*}
again by maps induced by inclusion.  The isomorphisms also hold for
the trivial case $k=n$.  Putting all this together yields
\begin{multline*}
\dsum_{k=0}^n 
    H(\dpair Qkn)[-2k]
\\\isom
\dsum_{u_\in Q'_{[0,n]}}
    H(\Delta[u,V],\Delta[u,V)\unite\Delta(u,V])[-2d(u)]
\end{multline*}
as a description of the cohomology of the complement of the affine
arrangement in $PV\wo PA_0$.
\end{rem}

The following wil be needed when deriving product information from the
spectral sequence in \prettyref{sec:graded-ring}.

\begin{prop} \label{prop:X-filtration}
Let $F$ be the filtration on the total complex $X$ 
defined by $(F_s X)_t\deq\dsum_{q=0}^s X_{t-q,q}$.
Then 
\begin{multline*}
\im \left(H_i(F_s X)\xto{\eps_\ast}H_i(PV,\Union P\A)\right)
    =\\= \dsum_{2k\ge i-s}
    \im\left(H_{i-2k} (\dpair Qkn)\xto{h_k}H_i(PV,\Union P\A)\right).
\end{multline*}
\end{prop}

\begin{proof}
Setting $(F_s Y)_t\deq\dsum_{q=0}^s Y_{t-q,q}$, the proof of
\prettyref{prop:X-Y} shows that $g$ also induces isomorphisms
$g_\ast\colon H_\ast(F_s Y)\xto\isom H_\ast(F_s X)$.  
Therefore 
\begin{align*}
\begin{split}
&\im \left(H_i(F_s X)\xto{\eps_\ast}H_i(PV,\Union P\A)\right)
       = \im \left(H_i(F_s Y)\xto{\eps_\ast\cmps g_\ast}
                H_i(PV,\Union P\A)
          \right)
    \\&\qquad
       = \im \left(H_i\left(\dsum_{q=0}^s Y_{\bullet-q,q}
                     \right)\xto{(\eps\cmps g)_\ast}
                H_i(PV,\Union P\A)
          \right)
    \\&\qquad
       = \im\left(
            \dsum_{k\colon 0\le i-2k\le s} H_{i-2k}(\dpair Qkn)
	    \xto{\sum_k h_k}
	    H_i(PV,\Union P\A)
          \right)
\end{split}
\end{align*}
as claimed.
\end{proof}

\section{The graded ring structure}
\label{sec:graded-ring}
In this section we show the following graded version
of~\prettyref{eqn:intro-product}.
\begin{prop}    \label{prop:graded-product}
Let $c\in H_\ast(Q_{[k,n]},Q_{[k,n)})$, 
$d\in H_\ast(Q_{[l,n]},Q_{[l,n)})$, $k+l\ge n$.
Then
\begin{equation*}
 h_k(c)\iprod h_l(d)=h_{k+l-n}(c\mintimes d)
 +\sum_{i\colon k+l-n<i\le n} h_i(r_i)
\end{equation*}
for classes $r_i\in H_\ast(Q_{[i,n]},Q_{[i,n)})$.
\end{prop}

To prove this proposition we will make the spectral sequence of the
preceding section into a multiplicative spectral sequence.  We will
use simplicial cohomology to achieve this.

Let $M$ be a triangulation of $PV$, which is a barycentric subdivision
of a triangulation of which all $A\in P\A$ are subcomplexes.  In
particular, all $Pu$, $u\in Q$, are full subcomplexes of $M$.
We will denote
the face poset of this triangulation by $FM$ and by $C(FM)$ the chain
complex of ascending (from $0$-simplices to $2n$-simplices) chains in
$FM$.  
\begin{defprop}
For a subcomplex $A$ of~$M$, cap
products
\begin{align*}
	C^r(FM, FM\wo FA)\tens C_s(FM)&\xto{\cap} C_{s-r}(FA) \\
        C^r(FM\wo FA)\tens C_s(FM)&\xto\cap C_{s-r}(FM, FA)\\
	h\tens{\simp{f_0,\dots,f_s}}&\mapsto 
	   (-1)^{r(s-r)}h(\simp{f_{s-r},\dots,f_s})\simp{f_0,\dots,f_{s-r}}
\end{align*}
are defined, where
\begin{align*}
C^r(FM, FM\wo FA)
  &=\ker
  \left(\Hom(C_r(FM), R)\to \Hom(C_r(FM\wo FA), R)
  \right),\\
C_{s-r}(FM, FA)
  &=\coker
  \left(C_{s-r}(FA)\to C_{s-r}(FM)\right).
\end{align*}
\end{defprop}
\begin{proof}
If $\simp{f_0,\dots,f_{s-r}}$ is not in $C(FA)$, 
then $f_{s-r}$ is not in $A$
and therefore $\simp{f_{s-r},\dots,f_s}$ is in $C(FM\wo FA)$, so
that $h(\simp{f_{s-r},\dots,f_s})=0$ for the first kind of product or
$h(\simp{f_{s-r},\dots,f_s})$ is defined for the second kind of product.
\end{proof}\sloppy
Since $\Delta(FM)$ is just the barycentric subdivison of~$M$,
we have $H(C(FM))\isom H(M)=H(PV)$.  Let $o\in C_{2n}(FM)$ represent the
orientation class $[M]\in H_{2n}(M)$.  
Regarding $C(FA)$ as a subcomplex of the singular complex~$S(A)$, this
yields a map $C(FM, FM\wo FA)\xto{\cap o} S(A)$ which induces an
isomorphism in homology, if $A$ is a full
subcomplex.  
$\Delta(FM\wo
FA)$ is the subcomplex of the barycentric subdivision of~$M$ that
consists of all simplices which do not meet~$A$.  It is the complement
of an open normal neighbourhood of~$A$, if $A$ is a full subcomplex.
\begin{defprop}
If $A$, $B$ are
subcomplexes of $M$, a cup product
\begin{align*}
	C^r(FM, FM\wo FA)\tens C^s(FM, FM\wo FB)
		&\to C^{r+s}(FM, FM\wo F(A\intersect B))\\
	g\tens h&\mapsto g\cup h,
\end{align*}
$(g\cup h)(\simp{f_0,\dots,f_{r+s}})\deq
(-1)^{rs}g(\simp{f_0,\dots,f_r})h(\simp{f_r,\dots,f_{r+s}})$,
is defined.
\end{defprop}
\begin{proof}
If $\simp{f_0,\dots,f_{r+s}}$ is in $C(FM\wo F(A\intersect B))$, then
$f_0\notin F(A\intersect B)$ and therefore either $f_0\notin FA$ and
$\simp{f_0,\dots,f_r}\in C(FM\wo FA)$, or $f_0\notin FB$ and 
$\simp{f_r,\dots,f_{r+s}}\in C(FM\wo FB)$. In either case 
$(g\cup h)(\simp{f_0,\dots,f_{r+s}})=0$. 
\end{proof}

We now introduce the double complex that will serve as a dual of the
double complex~$X$.

\begin{defn}
We define a double complex $Z$ by
\begin{equation*}
Z_{p,q}\deq\dsum_{\substack{u\in Q^{q+1}\\u_0<\dots<u_q=V}} 
    C^{-p}(FM, FM\wo F(Pu_0))
\end{equation*}
and, using notation as in \prettyref{def:X}, 
\begin{align*}
    \dd'\colon Z_{p,q}&\to Z_{p-1,q}\\
    f\tens\simp{u_0,\dots,u_q}&\mapsto\delta f\tens\simp{u_0,\dots,u_q}\\
\intertext{and}
    \dd''\colon Z_{p,q}&\to Z_{p,q-1}\\
    f\tens\simp{u_0,\dots,u_q}&
    \mapsto\sum_{i=0}^{q-1}(-1)^{p+i}f
    \tens\simp{u_0,\dots,\hat{u_i},\dots,u_q},
\end{align*}
satisfying $\dd'\dd''+\dd''\dd'=0$.
\end{defn}

\begin{defn}
We define products on $Z$ by
\begin{align*}
  Z_{p,q}\tens Z_{p',q'} &\to Z_{p+p',q+q'} \\
  (f\tens c)
    \tens(g\tens d)
  &\mapsto
    (-1)^{qp'}
    (f\cup g)\tens(c\mintimes d)
\intertext{and on $Y$ by}
  Y_{p,q}\tens Y_{p',q'} &\to Y_{p+p'-n,q+q'} \\
  c\tens d &\mapsto c\mintimes d.
\end{align*}
\end{defn}

\begin{defprop} \label{prop:Z-mult}
The map
\begin{align*}
g'\colon Z_{p,q}&\to X_{p+2n,q} \\
f\tens c&\mapsto (f\cap o)\tens c
\end{align*}
leads to a commutative diagram
\begin{equation*}\xymatrix@C+2em{
    H_\ast(Z)
    \ar[r]_\isom^{g'_\ast}
    \ar[d]
    &
    H_\ast(X)
    \ar[d]^{\eps_\ast}_\isom
    \\
    H^\ast(PV\wo\Union P\A)
    \ar[r]^-{\cap[PV]}_\isom
    &
    H_\ast(PV, \Union P\A),
}\end{equation*}
where the arrow on the left hand side is a ring homomorphism.
\end{defprop}

\begin{proof}
The diagram
\begin{equation*}\xymatrix@C+1em{
    Z_{p,q}
    \ar[r]^{g'}\ar[d]^{\eps'}
    &
    X_{p+2n,q}
    \ar[d]^\eps
    \\
    C^{-p-q}(FM\wo F(\Union P\A))
    \ar[r]^-{\cap o}
    &
    S_{2n+p+q}(PV,\Union P\A),
    \\
    S^{-p-q}(PV\wo\Union P\A)
    \ar[u]\ar[ru]_{\cap o}
}\end{equation*}
where
\begin{align*}
    \eps'\colon Z_{p,q} 
        &\to C^{-p-q}(FM\wo F(\Union P\A))\\
    f\tens\simp{u_0,\dots,u_q}
        &\mapsto0,\qquad\text{if $q>0$,}\\
    f\tens\simp V
        &\mapsto f,
\end{align*}
is commutative.
When checking that $\eps'$ is a chain map from the total complex, the
interesting case is $Z_{p,1}$.  Here
\begin{equation*}
\eps'(\dd(f\tens\simp{u,V}))
    = \eps'(\delta f\tens\simp{u,V}+(-1)^p f\tens\simp V)
    = (-1)^p f,
\end{equation*}
and since $f\in C^{-p}(FM, FM\wo F(Pu_0))$, the restriction of~$f$ to 
the chains from $C(FM\wo F(\Union P\A))$ is zero.  
That $\eps'$ respects products
is obvious.

The map $g'$ induces an isomorphism 
$\Hp_\ast(Z_{\bullet,q})\to\Hp_\ast(X_{\bullet,q})$ 
and therefore an
isomorphism $H_\ast(Z)\to H_\ast(X)$ between the total complexes.  The map 
$H^\ast(PV\wo\Union P\A)\to H^\ast(C(FM\wo F(\Union P\A)))$
induced by the lower left arrow is also an 
isomorphism.
\end{proof}

\begin{prop} \label{prop:Z-Y-mult}
The two isomorphisms
\begin{equation*}\xymatrix{
  E^1(Z)\ar[r]_-\isom^-{g'_\ast}
  &
  E^1(X)
  &
  E^1(Y)\ar[l]^-\isom_-{g_\ast}
}\end{equation*}
induce the same multiplication on $E^1(X)$.
\end{prop}

\begin{proof}
As noted in the proof of \prettyref{prop:X-Y}, we have 
\[E^1_{p,q}(X)\isom\Hp_p(X_{\bullet, q})\isom
    \dsum_{u_0<\dots<u_q=V} H_p(Pu_0)
\]
and
\[
  H_p(Pu_0) \isom 
  \begin{cases}
    \Z,&\text{$p=2k$, $0\le k\le d(u_0)$,}\\
    0,&\text{otherwise}.
  \end{cases}
\]
For $u\in Q_{[k,n]}$, there is a well defined generator 
$e^u_k\in H_{2k}(Pu)$ with the property
that $e^u_k$ is the image of the
canonical orientation class $[\C P^k]$ under a complex embedding 
$\C P^k\to Pu$.  In particular the map
\[
Y_{2k,q}=\Hp_{2k}(Y_{\bullet, q})\isom E^1_{2k,q}(Y)
  \xto[\isom]{g_\ast}
  E^1_{2k,q}(X)\isom \dsum_{u_0<\dots<u_q=V} H_{2k}(Pu_0)
\]
takes $\simp{u_0,\dots,u_q}$ to 
$e^{u_0}_k\tens\simp{u_0,\dots,u_q}$ and induces the multiplication
\begin{multline*}
  (e^{u_0}_k\tens\simp{u_0,\dots,u_q})
  \tens
  (e^{v_0}_l\tens\simp{v_0,\dots,v_{q'}})
  \\\mapsto
  \begin{cases}
    e^{u_0\lmin v_0}_{k+l-n}\tens(\simp{u_0,\dots,u_q}\mintimes
                                  \simp{v_0,\dots,v_{q'}}),
    &k+l\ge n,\\
    0,&k+l<n
  \end{cases}
\end{multline*}
on $E^1(X)$.

On the other hand, the multiplication on $E^1(X)$ induced via $g'$ by that 
on~$Z$ is seen to be
\begin{equation*}
  (a\tens c)\tens(b\tens d)
  \mapsto (-1)^{\adeg c+\adeg b}
          (a\iprod b) \tens (c\mintimes d),
\end{equation*}
where the intersection product is defined by commutativity of
\begin{equation*}\xymatrix{
    H^\ast(PV,PV\wo Pu) \tens H^\ast(PV,PV\wo Pv)
    \ar[r]^-\cup
    \ar[d]^{(\cap [PV])\tens(\cap [PV])}_\isom
    &
    H^\ast(PV, PV\wo (Pu \intersect Pv))
    \ar[d]^{\cap[PV]}_\isom
    \\
    H_\ast(Pu)\tens H_\ast(Pv)
    \ar[r]^-\iprod
    &
    H_\ast(Pu\intersect Pv),
}\end{equation*}
and we have
\begin{equation*}
   e^u_k\iprod e^v_l =
   \begin{cases}
     e^{u\lmin v}_{k+l-n},& k+l\ge n,\\
     0,&k+l<n,
   \end{cases}
\end{equation*}
so that the two products agree.
\end{proof}

\begin{proof}[Proof of \prettyref{prop:graded-product}]
By \prettyref{prop:Z-mult} the multiplication on $H(X)$ 
induced by that on $Z$
is carried by $\eps_\ast$ into the intersection product on 
$H_\ast(PV,\Union P\A)$. By \prettyref{prop:Z-Y-mult} 
the complexes $Z$ and $Y$ have
multiplicative spectral sequences which are isomorphic from 
the $E^1$-terms on,
and these induce a multiplicative structure on the spectral sequence
of $X$.  The $E^\infty$-term of this spectral sequence (which
degenerates at the $E^2$-term) is the graded object corresponding to
the filtration considered in \prettyref{prop:X-filtration}.  By
\prettyref{prop:Z-Y-mult} the $E^\infty$-term of $X$ is isomorphic to
that of~$Y$.
\end{proof}

\section{Intersecting with a hyperplane}
\label{sec:hyperplane}
Let $\Lambda\colon V\to\C$ be a linear functional that vanishes on no
element of $Q_{[0,n]}$ and $H\deq\ker\Lambda$.  $\A$ induces an
arrangement $\A^H\deq\set{A\intersect H\colon A\in\A}$ in $H$.  
If we
denote the intersection poset of $\A^H$ by $Q^H$, 
\begin{align*}
\eta\colon Q_{(0,n]}&\to Q^H_{[0,n-1]}\\
   q&\mapsto q\intersect H
\end{align*}
is an isomorphism lowering dimensions by one.

\begin{prop}  \label{prop:hk-barhk}
Consider the inclusion map 
$i\colon (PH, \Union P\A^H)\to (PV,\Union P\A)$.
For $c\in H_\ast(\dpair Qkn)$
we have
\begin{equation*}
i_!(h_k(c))=\begin{cases}
                   h^H_{k-1}(\eta_\ast(c)),& k>0,\\
		   0,&k=0.
	    \end{cases}
\end{equation*}
In particular $\ker i_!=\im h_0$.
\end{prop}

\begin{proof}
We first choose $(x^u_j)_{0\le j\lt k}$ with $x^u_j\in u\intersect H$
satisfying the conditions of \prettyref{def:fk} and therefore defining
functions $f_H^{k-1}$ and $h^H_{k-1}$.  Now for each $u\in Q_{[k,n]}$
we choose $x^u_k\in u$ with $\Lambda(x^u_k)=1$.  $(x^u_j)_{0\le j\le
k}$ then also satisfies the conditions of \prettyref{def:fk} and can
be used to define $f^k$ and $h_k$.

Indeed we calculate
\begin{equation*}
    \Lambda\left(\sum_{j=0}^k\sum_{i=0}^r\mu_j\lambda_i x^{u_i}_j\right)
    =    \sum_{i=0}^r\mu_k\lambda_i
    =  \mu_k.
\end{equation*}
First this 
implies that $f^k(x,y)\in H$ iff $x\in\C P^{k-1}\subset \C P^k$.
In particular $f^0$ misses $H$, which proves that part of the
proposition, and we now assume $k>0$.
The equation also implies 
that for $x=[x_0:\dots:x_{k-1}]\in\C P^{k-1}$ and 
$y\in\Delta(Q_{[k,n]})$ the map 
$\mu\mapsto f^k([x_0:\dots:x_{k-1}:\mu],y)$ meets $H$ transversally.
Furthermore
\begin{align*}
f^k\left([\mu_0:\dots:\mu_{k-1}:0],\sum_{i=0}^r \lambda_i u_i\right)
	&=
		\left[\sum_{j=0}^{k-1}\sum_{i=0}^r\mu_j\lambda_i x^{u_i}_j
		\right]
        \\&= f_H^{k-1}\left([\mu_0:\dots:\mu_{k-1}], 
	               \sum_{i=0}^r \lambda_i \eta(u_i)\right),
\end{align*}
which proves the proposition as we will now show in more detail.

\def\thom{\vartheta}
Let $\thom\in H^2(PV,PV\wo PH)$ be the Thom class of (the normal
bundle of) $PH$ in $PV$, i.e.\ the class satisfying
$\theta\cap[PV]=[PH]$.  
By the above calculations 
$(f^k)^\ast(\thom)\in H^2\left(\C P^k\times\Delta Q_{[k,n]},
                          (\C P^k\wo\C P^{k-1})\times\Delta
			  Q_{[k,n]}\right)$
is the Thom class of $\C P^{k-1}\times\Delta Q_{[k,n]}$ in 
$\C P^k\times\Delta Q_{[k,n]}$ which equals the class 
$\alpha\times 1$ where $\alpha\in H^2(\C P^k,\C P^k\wo\C P^{k-1})$
which is again a Thom class and maps to the canonical generator of
$H^\ast(\C P^k)$.  We finally calculate
\begin{align*}
i_!(h_k(c))
&=\thom\cap h_k(c)
\\&=\thom\cap f^k_\ast([\C P^k]\times c)
\\&=f^k_\ast\left((f^k)^\ast(\thom)\cap([\C P^k]\times c)\right)
\\&=f^k_\ast\left((\alpha\times1)\cap([\C P^k]\times c)\right)
\\&=f^k_\ast\left((\alpha\cap[\C P^k])\times (1\cap c)\right)
\\&=f^k_\ast\left([\C P^{k-1}]\times c\right)
\\&=h_{k-1}(\eta_\ast(c))
\end{align*}
as claimed.
\end{proof}

This proposition will enable us to prove statements about arrangements
by induction on the dimension of~$V$.  An example is the derivation of
\prettyref{thm:main} from Propositions~\ref{prop:graded-product}
and~\ref{prop:vanishing} in \prettyref{sec:main-proof}.

\section{Real projective arrangements and an example}
\label{sec:real}
We will investigate, without giving full proofs, what of the above
remains true for real projective arrangements, see how the analogue of
the product formula we are trying to prove fails for real projective
arrangements, and sketch the difference between real and complex
arrangements that will allow us to prove the formula for complex
arrangements.

For a real projective arrangement, functions $f^k$ as in
\prettyref{def:fk} exist but are not uniquely defined up to homotopy.
A direct sum decomposition as in \prettyref{prop:proj-homology} will
not be achieved for integral coefficients, but working over the
ring~$\Z_2$, we obtain an isomorphism
\begin{equation*}
\sum_{k=0}^n h_k\colon\dsum_{k=0}^n 
             H_\ast(\Delta Q_{[k,n]}, \Delta Q_{[k,n)};\Z_2)[-k]
	\xto{\,\isom\,} H_\ast(PV, \Union P\A;\Z_2).
\end{equation*}
The proof of \prettyref{prop:graded-product} will require no
additional changes.

Let $k,l\ge0$, $n\deq k+l+1$.  We consider the following subspaces of 
$\R^{n+1}=\R^k\times\R^l\times\R^2$.
\begin{align*}
    u&\deq\R^k\times\set0\times(\R\cdot(0,1)), &
    \tilde u&\deq\set0\times\R^l\times(\R\cdot(1,1)), \\
    v&\deq\R^k\times\set0\times(\R\cdot(4,1)), &
    \tilde v&\deq\set0\times\R^l\times(\R\cdot(5,1)),
\end{align*}
In the arrangement $\check\A\deq\set{u,v,\tilde u,\tilde v}$ we will
find classes $c\in H(\dpair Qkn)$, $d\in H(\dpair Qln)$ with
$h_k(c)\iprod h_l(d)\ne0$, although $k+l<n$.

The combinatorial data of $\check\A$ are given by the intersection poset
\begin{equation*}\xymatrix@C=0pt@R=0.5ex{
    &&&&&V
    \ar@{-}[dlll]
    \ar@{-}[dl]
    \ar@{-}[dr]
    \ar@{-}[drrr]
    \\
    \check Q_{[0,n]}=
    &&u&&v&&\tilde u&&\tilde v
    \\
    &&&u\intersect v
    \ar@{-}[lu]
    \ar@{-}[ru]
    &&&&\tilde u\intersect\tilde v
    \ar@{-}[lu]
    \ar@{-}[ru]
}\end{equation*}
with $u\intersect v$ and $\tilde u\intersect \tilde v$ only present
for $k>0$ and $l>0$ respectively,
and the dimensions $d(u)=d(v)=k$, $d(u\intersect v)=k-1$, 
$d(\tilde u)=d(\tilde v)=l$, $d(\tilde u\intersect\tilde v)=l-1$,
$d(V)=n$.  In case of $k=l$, we can, if we want
to, avoid the intersections $u\intersect v$ and 
$\tilde u\intersect\tilde v$ by a small change of $u$ and~$\tilde u$ without substantially affecting the
calculations below.  This shows that, in contrast to the case of affine
    arrangements in
\cite{delgormacph}
and \cite{delong-schultz}, a simple condition on the occuring
codimensions will not be enough for the product formula to extend from
    complex to real arrangements.

To simplify the pictures below and to have the notation parallel that
of \prettyref{sec:vanishing}, we consider $\check A$ to be the
union of the two arrangements 
$\A\deq\set{u,v}$ and $\tilde A\deq\set{\tilde u,\tilde v}$.
The arrangement $\tilde\A$ is in
general position with respect to $\A$ as in
\prettyref{def:general-position}.

Denoting the intersection posets of $\A$ and~$\tilde\A$ by $Q$
and~$\tilde Q$ respectively, we have $H_1(\dpair Qkn;\Z_2)\isom\Z_2$,
generated by $c\deq [\simp{u,V}+\simp{v,V}]$, and $H_1(\dpair{\tilde
Q}ln;\Z_2)\isom\Z_2$, generated by $d\deq [\simp{\tilde
u,V}+\simp{\tilde v,V}]$.
For the definition of
\begin{align*}
f^k\colon\R P^k\times\dpairp Qkn&\to (V,\Union P\A),\\
\tilde f^l\colon\R P^l\times\dpairp {\tilde Q}kn&\to (V,\Union P\tilde\A),
\intertext{and hence of}
h_k\colon H_r(\dpair Qkn;\Z_2)&\to H_{r+k}(V,\Union P\A;\Z_2),\\
\tilde h_l\colon H_r(\dpair{\tilde Q}kn;\Z_2)
    &\to H_{r+l}(V,\Union P\tilde\A;\Z_2)
\end{align*}
we set, with $(e_0,\dots,e_n)$ the standard basis of $V=\R^{n+1}$,
\begin{align*}
x^u_j&\deq\begin{cases}
             e_j,&j<k,\\
	     e_{k+l+1},&j=k,
          \end{cases}
&
\tilde x^{\tilde u}_j&\deq\begin{cases}
             e_{k+j},&j<l,\\
	     e_{k+l}+e_{k+l+1},&j=l,
          \end{cases}
\\
x^v_j&\deq\begin{cases}
             e_j,&j<k,\\
	     4e_{k+l}+e_{k+l+1},&j=k,
          \end{cases}
&
\tilde x^{\tilde v}_j&\deq\begin{cases}
             e_{k+j},&j<l,\\
	     5e_{k+l}+e_{k+l+1},&j=l,
          \end{cases}
\\
x^V_j&\deq\begin{cases}
             e_j,&j<k,\\
	     2e_{k+l}+e_{k+l+1},&j=k,
          \end{cases}
&
\tilde x^{V}_j&\deq\begin{cases}
             e_{k+j},&j<l,\\
	     3e_{k+l}+e_{k+l+1},&j=l.
          \end{cases}
\end{align*}
To determine $h_k(c)\iprod h_l(d)$, we first have a look at the geometric
intersection 
$S\deq f^k\bigl[\R P^k\times\Delta Q_{[k,n]}\bigr]\intersect
    f^l\bigl[\R P^l\times\Delta \tilde Q_{[l,n]}\bigr]$.
For $x\in\Delta (Q_{[k,n]})$, $y\in\Delta (\tilde Q_{[l,n]})$, 
the intersection 
$f\left[\R P^k\times\set x\right]\intersect
 f\left[\R P^l\times\set y\right]$ is either empty or consists of a
single point.  The left of the following two pictures shows the two
dimensional simplicial complex 
$\Delta Q_{[k,n]}\times\Delta \tilde Q_{[l,n]}
 =\Delta(Q_{[k,n]}\times\tilde Q_{[l,n]})$.
\begin{equation}\label{eq:poset-drawing}
\begin{split}
\xymatrix@C=0.6em@R=0.9em{
  (u,\tilde v)
    \ar@{-}[rr]
    \ar@{-}[dd]
  &&(V,\tilde v)
  &\bullet&(v,\tilde v)
    \ar@{-}[ll]
    \ar@{-}[dd]
  &&
  &&\tilde v
  &&
  \\
  \\
  (u,V)
  &&(V,V)
    \ar@{-}[lluu]
    \ar@{-}[uu]
    \ar@{-}[rruu]
    \ar@{-}[rr]
    \ar@{-}[rrdd]
    \ar@{-}[dd]
    \ar@{-}[lldd]
    \ar@{-}[ll]
  &&(v,V)
  &&
  u\ar@{.}@<3pt>[rr]
  &&V
    \ar@{-}[uu]
    \ar@{-}[rr]
    \ar@{-}[dd]
    \ar@{-}[ll]
    \ar@{.}@<3pt>[uu]
  &&v
  \\
  \bullet\ar@{.}[rrruuu]
  \\
  (u,\tilde u)
    \ar@{-}[rr]
    \ar@{-}[uu]
  &&(V,\tilde u)
  &&(v,\tilde u)
    \ar@{-}[ll]
    \ar@{-}[uu]
  &&
  &&\tilde u
  &&
}
\end{split}
\end{equation}
The dotted line depicts the set $\bar S$ of those points $(x,y)$ 
for which the
intersection is nonempty.  $S$~is a connected $1$--dimensional
manifold with boundary, and we can see from the picture that one
boundary point lies in $u\intersect V=u$ and the other one in
$V\intersect\tilde v=\tilde v$.  A closer look at~$S$, which is the
intersection of two manifolds that meet transversely, shows that
indeed 
$h_k(c)\iprod \tilde h_l(d)=
    \check h_0([\simp{u,V}+\simp{\tilde v,V}])$.  
$[\simp{u,V}+\simp{\tilde v,V}]$ is a generator of
$H_1(\dpair{\check Q}0n;\Z_2)$, therefore 
$h_k(c)\iprod \tilde h_l(d)\ne0$.

We equip $Q\times\tilde Q$ with a dimension function 
$d(p,q)\deq d(p)+d(q)$.
The map 
$Q\times\tilde Q\to\check Q$, $(p,q)\mapsto p\intersect q$, sends
$\pairp{(Q\times\tilde Q)}n{2n}$ to $\pairp{\check Q}0n$.
In the picture, the border of the square is 
$\Delta(Q_{[k,n]}\times\tilde Q_{[l,n]})_{[0,2n)}$ and the four vertices at
the corners are $\Delta(Q_{[k,n]}\times\tilde Q_{[l,n]})_{[0,n)}$.  Under
the composition of maps
\begin{equation}\label{eq:oder-complex-retract}
\begin{split}\xymatrix@R-1ex{
\genfrac{}{}{0pt}0
 {\bigl(\Delta(Q_{[k,n]}\times\tilde Q_{[l,n]})\wo
    \Delta(Q_{[k,n]}\times\tilde Q_{[l,n]})_{[0,n)},
 }{\qquad\qquad\quad
 \Delta(Q_{[k,n]}\times\tilde Q_{[l,n]})_{[0,2n)}\wo
    \Delta(Q_{[k,n]}\times\tilde Q_{[l,n]})_{[0,n)}
 \bigr)
 }
\ar[d]\\
 \dpairp{(Q_{[k,n]}\times\tilde Q_{[l,n]})}n{2n}
\ar[d]\\
  \dpairp{\check Q}0n,
}\end{split}\end{equation}
the first being a deformation retraction and the second given by
inclusion, in our example $(\bar S,\dd\bar S)$ is mapped to the dotted
line in the picture on the right.  This set carries the relative cycle 
$\simp{u,V}+\simp{\tilde v,V}$ representing 
$h_k(c)\iprod \tilde h_l(d)$.  We will see in
 \prettyref{sec:vanishing} that this is not just a coincidence.

When considering complex arrangements we will see that in the above
situation we gain one dimension compared to real arrangements, and
$\bar S$ will miss the cone with top the vertix $(V,V)$ and base
$\Delta(Q_{[k,n]}\times\tilde Q_{[l,n]})_{[0,n)}$.  The map of 
$(\bar S, \dd\bar S)$ to $\dpairp{\check Q}0n$ 
will therefore miss the vertix $V$ and be homotopic to a map with
image in $\Delta\check Q_{[0,n)}$.  This is the idea behind the proof
of \prettyref{prop:h0-vanish}, although it will be technically a bit
different.

\section{Recovering the direct sum decomposition}
\label{sec:stratum-top}
When discussing the real example, it seemed plausible that a certain
subset of the order complex of the intersection poset \emph{should}
carry the inverse image of the considered intersection product under
the isomorphism $\sum_k h_k$.  We now develop tools that allow to
actually \emph{prove} this kind of proposition.

More generally, given a class in $H_\ast\left(PV, \Union P\A\right)$
we want to determine the corresponding element of $\dsum_k
H_\ast\dpairp Qkn$.  Because of \prettyref{prop:hk-barhk} it will
suffice to identify the part in the summand $H_\ast\dpairp Q0n$.  The
key to this will be to not only consider the map $f^0\colon \Delta
Q_{[0,n]}\to PV$, but also a map $PV\to Q_{[0,n]}$, where the poset
$Q$ is topologized in an appropriate way yielding the \emph{space of
strata}.  While we have up to this point used only the former map, in
\cite{delgormacph} a description of the cohomology ring of the
complement of an affine arrangement is obtained using exclusively the
latter map.  Here the interplay of both maps will be important.

\begin{defn}
Let $P$ be a poset.  We make $P$ into a topological space by calling a
set~$O\subset P$ open, iff $x\in O$ implies $y\in O$ for all $y\ge x$.
\end{defn}

\begin{lem}     \label{lem:poset-homot}
Let $X$ be a space, $P$ a poset, $A\subset X$, $R\subset P$.  If
$f,g\colon (X,A)\to(P,R)$ are continuous maps with $f(x)\ge g(x)$ for
all $x\in X$, then $f\homot g$.
\end{lem}

\begin{proof}
The desired homotopy is given by
\begin{align*}
H\colon (X,A)\times I&\to (P,R) \\
(x,t)&\mapsto
\begin{cases}
  f(x),& t<1,\\
  g(x),& t=1.
\end{cases}
\end{align*}
This map is continuous, since $g^{-1}[O]\subset f^{-1}[O]$ 
for open $O\subset X$, and therefore
$H^{-1}[O]=f^{-1}[O]\times[0,1)\unite g^{-1}[O]\times\set1
  = f^{-1}[O]\times [0,1)\unite g^{-1}[O]\times I.
$
\end{proof}

\begin{lem}
If $P$ has a minimum or a maximum, then $P$ is contractible.
\end{lem}

\begin{proof}
By the preceding lemma, the constant map to the minimum respectively
the maximum is homotopic to the identity.
\end{proof}

\begin{defn}
\def\C{{\mathfrak C}}
Let $X$ be a space equipped with a covering $\C$ by closed sets and
let $P$ be the poset $P\deq\set{\Intersection M\colon \emptyset\ne
M\subset\C,\Intersection M\ne\emptyset}$, odered by inclusion.  
We define a continous map
\begin{align*}
  s\colon X &\to P \\
  x & \mapsto \min\set{p\in P\colon x\in P}=\Intersection\set{C\in\C\colon
  x\in C}.
\end{align*}
In particular we consider the following two kinds of maps.
For our arrangement $\A$ we consider the map
$s^\A\colon PV\to Q_{[0,n]}$ corresponding to the covering
$P\A\unite\set{P V}$ of $PV$.
For a poset $P$ which has unique minima in the sense that for
$M\subset P$, $M\ne\emptyset$,
the set $\set{p\in P\colon \text{$p\le q $ for all $q\in
M$}}$ is either empty or of the form $\set{p\colon p\le q}$ for a
$q\in P$, we consider 
the map $s^P\colon
\Delta P\to P$ arising from the covering of $\Delta P$ by the
subspaces $\Delta(\set{p'\colon p'\le p})$, $p\in P$.
\end{defn}

\begin{lem}   \label{lem:DP-P}
For a finite poset $P$ and $R\subset P$, both satisfying the condition
regarding minima of the preceding definition, the map $s^P_\ast\colon
H_\ast(\Delta P, \Delta R)\xto\isom H_\ast(P, R)$ is an isomorphism.
\end{lem}

\begin{proof}
We may assume $R=\emptyset$, because the general case will follow by
an application of the five lemma.

If $P=\emptyset$, $H_\ast(\Delta P)\isom0\isom H_\ast(P)$.

If $P\ne\emptyset$, $P$ has a minimal element $m$.  We set
$M\deq\set{p\colon p\ge m}$.  $P=M\unite (P\wo\set m)$ is an open
covering, and $\Delta P = \Delta M \unite \Delta(P\wo\set m)$ a
covering by subcomplexes, and $\Delta M \intersect\Delta(P\wo\set
m)=\Delta(M\wo\set m)$.  By induction on the number of elements, the
maps $H_\ast(\Delta (P\wo\set m))\to H_\ast(P\wo\set m)$ and
$H_\ast(\Delta (M\wo\set m))\to H_\ast(M\wo\set m)$ are isomorphisms.
$H_\ast(\Delta M)\to H_\ast(M)$ is an isomorphism, because both spaces
are contractible.  It follows by the Mayer-Vietoris theorem and the
five-lemma, that $H_\ast(\Delta P)\to H_\ast(P)$ is also an isomorphism.
\end{proof}

\begin{prop} \label{prop:detect-h0}
The composition
\[
H_\ast(\dpair Qkn)
  \xto{h_k}H_\ast(PV,\Union P\A) \xto{s^\A_\ast}
  H_\ast(Q_{[0,n]},Q_{[0,n)})
\]
is an isomorphism for $k=0$ and zero for $k>0$.
\end{prop}

\begin{proof}
Consider the diagram
\begin{equation*}\xymatrix@C-0.5em{
    \C P^k\times\dpairp Qkn
    \ar[rr]^-{f^k}
    \ar[d]^\pi
    &
    &
    (PV,\Union P\A)
    \ar[d]^{s^\A}
    \\
    \dpairp Qkn
    \ar[r]^-{s^{Q_{[k,n]}}}
    &
    (Q_{[k,n]},Q_{[k,n)})
    \ar[r]^-i
    &
    (Q_{[0,n]},Q_{[0,n)})
}\end{equation*}
where $i$ is the inclusion map and $\pi$ the projection onto the
second factor.  By construction of $f^k$, 
$f^k\left[\set x\times\simp{u_0,\dots,u_q}\right]\subset Pu_q$, 
that is 
$s^\A(f^k(x,y))\le s^Q(y)$, and by \prettyref{lem:poset-homot} this
implies the homotopy commutativity of the diagram.

For $k=0$, $\pi$ is a homeomorhism, $i$ the identity, and $h_0$ equals
$f^0_\ast$ up to an isomorphism.  
Therefore $s^\A_\ast\cmps h_0$ is an isomorphism,
because $s^{Q_{[0,n]}}_\ast$ is an isomorphism by
\prettyref{lem:DP-P}.

For $k>0$, 
$s^\A_\ast(h_k(c))=s^\A_\ast(f^k_\ast([\C P^k]\times c))
    =(i\cmps s^{Q_{[k,n]}})_\ast(\pi_\ast([\C P^k]\times c))
    =(i\cmps s^{Q_{[k,n]}})_\ast(0)=0.$
\end{proof}

\section{The product of two arrangements}
\label{sec:arr-product}
The intersection of two sets can be identified with the intersection
of their cartesian product with a diagonal.  Similarly the
intersection product of two homology classes equals the image of their
cross product under the image of the transfer map associated with the
diagonal map.  We therefore study the products of two arrangements.

Already in the real example we have discussed, it was useful to assume
the homology classes of which the product was to be determined to be
carried by different arrangements. 
So we now assume to be given a second arrangement $\tilde A$ in $V$ with
intersection poset~$\tilde Q$.
For the first part of this section the arrangement could be in
a vector space different from~$V$, but we will have no use for this
generality later on.

We equip the poset 
$Q\times \tilde Q$ with a dimension function~$d$ by 
$d(u,v)\deq d(u)+d(v)$.
The counterpart of~$h$ for~$\tilde A$ will be denoted by~$\tilde h$
and so on.

As noted above, we will be interested in cross products.
\begin{defprop} \label{def:prod-g}
Any choice of $(y_i^{u,v})_{i=0,\dots,k}$, $(z_i^{u,v})_{i=0,\dots,l}$
for $(u,v)\in Q_{[k,n]}\times \tilde Q_{[l,n]}$ with $y^{u,v}_i\in u$,
$z^{u,v}_i\in v$ and such that for all
$(u_0,v_0)<\dots<(u_m,v_m)$ and $\lambda\in\Delta^m$ the system 
$(\sum_j\lambda_jy^{u_j,v_j}_i)_i$ as well as the system
$(\sum_j\lambda_jz^{u_j,v_j}_i)_i$ is linearly independent, yields a
map
\begin{align*}\scriptstyle
  g\colon
  \C P^k\times \C P^l\times
  \left(
    \Delta(Q_{[k,n]}\times\tilde Q_{[l,n]}),
    \Delta(Q_{[k,n]}\times\tilde Q_{[l,n]})_{[0,2n)}
  \right)
   &\scriptstyle
   \to (PV, \Union P\A)\times(PV,\Union P\tilde\A)\\
  \scriptstyle
  \left(
   [\mu_0:\dots:\mu_k],
   [\nu_0:\dots:\nu_l],
   \sum\limits_j\lambda_j (u_j,v_j)
  \right)
  &\scriptstyle\mapsto
  \left(
   \left[\sum\limits_{i,j}\lambda_j\mu_i y^{u_j,v_j}_i\right],
   \left[\sum\limits_{i,j}\lambda_j\nu_i z^{u_j,v_j}_i\right]
  \right),
\end{align*}
and as in \prettyref{def:fk}, any two such maps are homotopic.
\qed
\end{defprop}

\begin{prop} \label{prop:prod-g}
For $c\in H_\ast(\dpair Qkn)$, $d\in H_\ast(\dpair{\tilde Q}ln)$,
we have
$h_k(c)\times \tilde h_l(d)
 =g_\ast\left([\C P^k]\times[\C P^l]\times(c\times d)\right)$.
\end{prop}

\begin{proof}
For the choice $y^{u,v}_i=x^u_i$, $z^{u,v}_i=\tilde x^v_i$, we just get
the map $f^k\times\tilde f^l$ up to identification of 
$  \C P^k\times \C P^l\times
  \left(
    \Delta(Q_{[k,n]}\times\tilde Q_{[l,n]}),
    \Delta(Q_{[k,n]}\times\tilde Q_{[l,n]})_{[0,2n)}
  \right)
$ with
$ \C P^k\times(\dpair Qkn)
  \times
  \C P^l\times(\dpair{\tilde Q}ln)$.
\end{proof}

As noted after discussing the real example, it will be important to
control the codimension of a set corresponding to the dotted line in
\prettyref{eq:poset-drawing}.  We will now work towards this and start
with an algebraic lemma.
\begin{lem} \label{lem:codim}
Let $u$, $v$ be subspaces of $V$ in general position with respect to
each other, 
$\dim u=r\ge k+1$, $\dim v=s\ge l+1$, $\dim V=n+1$, $k+l<n$.
Let $O$ be the open subspace of the affine space 
$u^{k+1}\times v^{l+1}$ defined by 
\begin{multline*}
  O\deq\bigl\{(y_0,\dots,y_k,z_0,\dots,z_l)\colon 
  \\\text{$\dim(\spann\set{y_i})=k+1$,
        $\dim(\spann\set{z_i})=l+1$}\bigr\}
\end{multline*}
and algebraic subsets $\cdots\subset S_1\subset S_0\subset O$ defined by
\begin{equation*}
S_m\deq\set{(y_0,\dots,y_k,z_0,\dots,z_l)\colon
  \dim(\spann(\set{y_i}\unite\set{z_j})) < k+l+2-m}
\end{equation*}
Then $S_m\wo S_{m+1}$ is a complex submanifold of codimension
$(1+m)(n-k-l+m)$.
\end{lem}

\begin{proof}
We consider $(y_0,\dots,y_k,z_0,\dots,z_l)\in S_m\wo S_{m+1}$.
This implies $u+v=V$.  We set $Y\deq\spann\set{y_i}$, 
$t\deq n-k-s+\dim(Y\intersect v)$, and choose a basis 
$(e_0,\dots,e_n)$ of~$V$ such that 
$\spann\set{e_0,\dots,e_{r-1}}=u$,
$\spann\set{e_{n-s+1},\dots,e_n}=v$,
$\spann\set{e_t,\dots,e_{k+1}}=Y$.  Let $A$ be the 
$(n+1)\times(k+l+2)$--matrix with columns
$(y_0,\dots,y_k,z_0,\dots,z_l)$ expressed using this basis.
Elements of $O$ are represented by
matrices $A'=(a'_{ij})$ with $a'_{ij}=0$ for $r\le
i\le n$, $0\le j\le k$ and for $0\le i\le n-s$, $k+1\le j\le k+l+1$
such that the first $k+1$ and the last $l+1$ columns are linearly
independent, and $A=(a_{ij})$ has the additional property that the
first $t$~rows are zero.

There are sets
$I$ and $J$ with $\set{t,\dots,k+t}\subset I\subset\set{t,\dots,n}$, 
$\set{0,\dots,k}\subset J\subset\set{0,\dots,k+l+1}$ 
and $\card I=\card J=k+l+1-m$ such that
the matrix $B\deq(a_{ij})_{\text{$i\in I$,$j\in J$}}$ is regular.  
Similarly,
there exist $I'$, $J'$ with $I'\subset\set{t,\dots,n}\wo M$,
$\set{k+1,\dots,k+l+1}\subset J'\subset\set{0,\dots,k+l+1}$ and $\card
I'=\card J'=k+l+1-m$ such that the matrix $C\deq(a_{ij})_{\text{$i\in
I'$,$j\in J'$}}$ is regular.  

Let $U\subset O$ be a neighbourhood of $A$ such that for every
$A'=(a'_{ij})\in U$ the matrices $(a'_{ij})_{\text{$i\in I$,$j\in J$}}$
and $(a'_{ij})_{\text{$i\in I'$,$j\in J'$}}$ are regular.
Then an $A'\in U$ is in $S_m$
if and only if the equations
\begin{align*}
f_{i_0j_0}(A')&\deq
\det(a'_{ij})_{
  \substack{i\in I\unite\set{i_0}\\
  j\in J\unite\set{j_0}}}
=0
&\text{for all\ }&
	i_0\in I_0\deq\set{n+1-s,\dots,n}\wo I,
\\&&&
	j_0\in J_0\deq\set{0,\dots,k+l+1}\wo J
\\
\intertext{and}
g_{i_0j_0}(A')&\deq\det(a'_{ij})_{
  \substack{i\in I'\unite\set{i_0}\\
            j\in J'\unite\set{j_0}}}
=0
&\text{for all\ }&
	i_0\in \set{0,\dots,t-1},
\\&&&
	j_0\in J'_0\deq\set{0,\dots,k+l+1}\wo J'
\end{align*}
hold.  To see this, assume $A'\notin S_m$, i.e.\ $\rank A'>k+l+1-m$.  
If the rank of the matrix $A'$ with the first $t$ rows deleted is
greater than $k+l+1-m$, one of the functions $f_{i_0j_0}$ becomes
non-zero, otherwise one of the functions~$g_{i_0j_0}$.

Finally we compute
for $(i_1,j_1)\in I_0\times J_0$
\begin{align*}
\bigabs{\frac{\partial f_{i_0j_0}(A)}{\partial a_{i_1j_1}}}
  &=\begin{cases}
	\abs{\det B},& (i_0,j_0)=(i_1,j_1),\\
	0,&(i_0,j_0)\ne(i_1,j_1),
    \end{cases}
&(i_0,j_0)&\in I_0\times J_0,
\\
\frac{\partial g_{i_0j_0}(A)}{\partial a_{i_1j_1}}
  &=0,
&(i_0,j_0)&\in \set{0,\dots,t-1}\times J'_0
\intertext{and for $(i_1,j_1)\in \set{0,\dots,t-1}\times J'_0$}
\frac{\partial f_{i_0j_0}(A)}{\partial a_{i_1j_1}}
  &=0,
&(i_0,j_0)&\in I_0\times J_0,
\\
\bigabs{\frac{\partial g_{i_0j_0}(A)}{\partial a_{i_1j_1}}}
  &=\begin{cases}
	\abs{\det C},& (i_0,j_0)=(i_1,j_1),\\
	0,&(i_0,j_0)\ne(i_1,j_1),
    \end{cases}
&(i_0,j_0)&\in \set{0,\dots,t-1}\times J'_0
\end{align*}
and $\card{I_0\times J_0\unite \set{0,\dots,t-1}\times J'_0}
	=(n+1-\card I)(m+1)=(n-k-l+m)\cdot(m+1).$
\end{proof}

\begin{defn} \label{def:general-position}
We say that the arrangement $\tilde\A$ is in general position with
respect to the arrangement $\A$, if for all $u\in Q$
and $v\in\tilde Q$, we have $u\intersect v=0$ whenever $d(u)+d(v)<n$
and $d(u\intersect v)=d(u)+d(v)-n$ otherwise.  
\end{defn}

\begin{prop} \label{prop:singular-codim}
Let $k+l<n$, $D\subset PV\times PV$ be the diagonal and
$S\subset\Delta(Q_{[k,n]}\times Q_{[l,n]})$ be defined as the set of
all points~$x$ such that
$g\left[\C P^k\times\C P^l\times\set x\right]
\intersect D\ne\emptyset$.
For a generic choice of the points
$y^{u,v}_i$ and $z^{u,v}_i$ defining~$g$, the set $S$ intersects every
open simplex of $\Delta(Q_{[k,n]}\times Q_{[l,n]})$ in an algebraic
set of real codimension $2(n-k-l)$.
\end{prop}

\begin{proof}
In regard of \prettyref{lem:codim} all that is required is that for
each chain $(u_0,v_0)<\dots<(u_t,v_t)$ the affine plane in
$u_t^{k+1}\times v_t^{k+1}$ spanned by the $t+1$ points
$(y^{u_0,v_0},z^{u_0,v_0}),\dots,(y^{u_t,v_t},z^{u_t,v_t})$ meets the
algebraic set $S_0$ transversely.  Assuming that the affine plane
spanned by the first~$t$ of these points already meets $S_0$
transversely, this will be fulfilled for a generic choice of
$(y^{u_t,v_t},z^{u_t,v_t})\in u_t^{k+1}\times v_t^{k+1}$.
\end{proof}

\section{The vanishing of the intersection product of classes of
  degrees not adding up to~\texorpdfstring{$n$}n.}
\label{sec:vanishing}
In this section we will prove the following proposition.  It is one
half of~\prettyref{eqn:intro-product}, and the other half will be
derived from it using an inductive argument relying also on
\prettyref{prop:hk-barhk} and \prettyref{prop:graded-product}.
\begin{prop} \label{prop:vanishing}
  Let $c\in H_\ast(\dpair Qkn)$, 
  $d\in H_\ast(\dpair Qln)$, $k+l<n$.
Then
\begin{equation*}
 h_k(c)\iprod h_l(d)=0.
\end{equation*}
\end{prop}

We will prove the proposition in
three steps.

We would like to have the classes $h_k(c)$ and $h_l(d)$ represented by
chains as much as possible in general position with respect to each
other.  To this end we consider an arrangement that is the union of
two arrangements $\A$ and~$\tilde\A$ with intersection posets $Q$
and~$\tilde Q$ such that $\tilde\A$ is in general position with
respect to~$\A$ (see \prettyref{def:general-position}).

We will denote the intersection poset of the arrangement 
$\check A\deq\A\unite\tilde A$ by~$\check Q$ and so on.  
The map
\begin{equation}\label{eq:sigma-ntc}
\begin{split}
    \sigma\colon(Q\times\tilde Q)_{[n,2n]}
        &\to \check Q_{[0,n]} \\
    (u,v) &\mapsto u\intersect v
\end{split}
\end{equation}
is an isomorphism.

\begin{prop} \label{prop:h0-vanish}
  In the above situation, let $c\in H_\ast(\dpair Qkn)$, 
  $d\in H_\ast(\dpair{\tilde Q}ln)$, $k+l<n$.
Then
\begin{equation*}
 h_k(c)\iprod \tilde h_l(d)=
  \sum_{i>0} \check h_i(r_i)
\end{equation*}
for classes $r_i\in H_\ast\left(\dpair{\check Q}in\right)$.
\end{prop}

\begin{proof}
By \prettyref{prop:detect-h0} we will have to show
$s^{\check A}_\ast(h_k(c)\iprod \tilde h_l(d))=0$.  It will be in
doing so that we employ the ideas laid out in the discussion of the
real example.

We set $(X,A)\deq \C P^k\times \C P^l\times
  \left(
    \Delta(Q_{[k,n]}\times \tilde Q_{[l,n]}),
    \Delta(Q_{[k,n]}\times \tilde Q_{[l,n]})_{[0,2n)}
  \right)
$, 
$D\deq\set{(x,x)\in PV\times PV}$, 
$\compl D\deq (PV\times PV)\wo D$ and use the map~$g$
from \prettyref{def:prod-g}.
We denote the diagonal map $PV\to PV\times PV$ by~$\Delta$ and define
$\bar g\colon (g^{-1}[D], g^{-1}[D]\intersect A)\to (PV, \Union
P\A\unite\Union P\tilde\A)$ by $\Delta\cmps\bar g=g$.
Note that in the real example the projection of $g^{-1}[D]$ 
to the order complex is the set represented by a dotted line in
  \prettyref{eq:poset-drawing}.
We will first show $h_k(c)\iprod \tilde h_l(d)\in\im \bar g_\ast$ and
  then $s^{\check A}_\ast\cmps \bar g_\ast=0$.

There is a commutative
  diagram
\begin{equation*}\xymatrix@C=0pt{
    H^\ast(PV\times PV, \compl D)
    \tens
    H_\ast(X, A)
    \ar[r]^-{\id\tens g_\ast}
    \ar[d]^{g^\ast\tens\id}
    &
    \relax\genfrac{}{}{0pt}0{
    H^\ast(PV\times PV, \compl D)
    \tens}{
    H_\ast\left((PV,\Union P\A)\times (PV,\Union P\tilde\A)\right)
    }
    \ar[dd]^\cap
    \\
    H^\ast(X, X\wo g^{-1}[D])
    \tens
    H_\ast(X,A)
    \ar[d]^\cap
    \\
    H_\ast(g^{-1}[D], g^{-1}[D]\intersect A)
    \ar[r]^-{g_\ast}
    \ar[rd]_-{\bar g_\ast}
    &
    H_\ast\bigl(
           D, D\intersect(\Union P\A\times PV
                          \unite
			  PV \times \Union P\tilde\A)
          \bigr)
    \\
    &
    H_\ast(PV, \Union P\A\unite\Union P\tilde\A).
    \ar[u]^\isom_{\Delta_\ast}
}\end{equation*}
Regarding the existence of the cap products in this diagram and
    commutativity, note that we are entirely dealing with algebraic
    sets and polynomial maps.
Now, if $\vartheta\in H^\ast(PV\times PV, \compl D)$ is the Thom class
determined by $\vartheta\cap[PV\times PV]=\Delta_\ast([PV])$, then
\begin{align*}
h_k(c)\iprod \tilde h_l(d)
  &= \Delta_!(h_k(c)\times \tilde h_l(d))
\\&= \Delta_\ast^{-1}\left(\vartheta\cap(h_k(c)\times \tilde h_l(d))
     \right)
\\&= \Delta_\ast^{-1}\left(\vartheta\cap 
    g_\ast\bigl([\C P^k]\times[\C P^l]\times(c\times d)\bigr)
    \right)
\\&= \Delta_\ast^{-1}\left(g_\ast\left(g^\ast(\vartheta)\cap
    \bigl([\C P^k]\times[\C P^l]\times(c\times d)\bigr)
    \right)
    \right)
\\&= \bar g_\ast\left(g^\ast(\vartheta)\cap
    \bigl([\C P^k]\times[\C P^l]\times(c\times d)\bigr)
    \right).
\end{align*}
By construction of $g$, 
$\bar g\left(x,y,\sum_{j=0}^r\lambda_j (u_j,v_j)\right)
  \in u_r\intersect v_r$.
Firstly this implies that $g^{-1}[D]$ misses 
$\C P^k\times\C P^l\times\Delta(Q\times\tilde Q)_{[0,n)}$, 
and secondly from the reformulation
$s^{\check\A}(\bar g(x,y,z))
\le \sigma\left(s^{Q\times\tilde Q}(z)\right)$,
where $\sigma$ is the isomorphism from \prettyref{eq:sigma-ntc},
and \prettyref{lem:poset-homot} it can be seen
that the diagram
\begin{equation*}\def\objectstyle{\scriptstyle}\xymatrix@C=-7em{
    (g^{-1}[D], g^{-1}[D]\intersect A)
    \ar[ddr]_-\pi
    \ar[rr]^-{s^{\check\A}\cmps \bar g}
    &&(\check Q_{[0,n]}, \check Q_{[0,n)})
    \\&&
    \relax\left(
    (Q\times\tilde Q)_{[n,2n]},
    (Q\times\tilde Q)_{[n,2n)}
    \right)
    \ar[u]_\sigma
    \\
    &
    \relax\left(
    \Delta(Q_{[k,n]}\times\tilde Q_{[l,n]})\wo
        \Delta(Q\times\tilde Q)_{[0,n)},
    \Delta(Q_{[k,n]}\times\tilde Q_{[l,n]})_{[0,2n)}\wo
        \Delta(Q\times\tilde Q)_{[0,n)}
    \right)
    \ar[ru]_{s^{Q\times\tilde Q}},
}\end{equation*}
where $\pi$ denotes projection onto the third factor,
is homotopy commutative.  The two arrows on the right hand side of the
diagram should be compared to \prettyref{eq:oder-complex-retract}.

By \prettyref{prop:singular-codim} and because 
the subcomplex 
$\Delta\bigl(
  (Q_{[k,n]}\times\tilde Q_{[l,n]})\intersect (Q\times\tilde Q)_{[0,n)}
\bigr)$ 
has dimension at most
$n-1-k-l$, we may assume that
$\pi[g^{-1}[D]]$
will not only miss this subcomplex, but any
cone over it.  Therefore $\pi$ factorizes over the pair
\begin{multline*}
\bigl(
    \Delta(Q_{[k,n]}\times\tilde Q_{[l,n]})\wo
        \Delta((Q\times\tilde Q)_{[0,n)}\unite\set{(V,V)}),
\\
    \Delta(Q_{[k,n]}\times\tilde Q_{[l,n]})_{[0,2n)}\wo
        \Delta(Q\times\tilde Q)_{[0,n)}
    \bigr).
\end{multline*}
This pair is homeomorphic to 
$(\Delta(Q_{[k,n]}\times\tilde Q_{[l,n]})_{[0,2n)}\wo
        \Delta(Q\times\tilde Q)_{[0,n)})
 \times ([0,1),\set0)$ and has trivial homology.  
So 
$(s^{\check\A}\cmps \bar g)_\ast
=s^{Q\times\tilde Q}_\ast\cmps \pi_\ast=s_\ast\cmps 0=0$.
\end{proof}

\begin{prop} \label{prop:vanishing-gp}
  In the above situation
  let $c\in H_\ast(\dpair Qkn)$, 
  $d\in H_\ast(\dpair{\tilde Q}ln)$, $k+l<n$.
Then
\begin{equation*}
 h_k(c)\iprod \tilde h_l(d)=0.
\end{equation*}
\end{prop}

\begin{proof}
We choose a hyperplane~$H$ in~$V$ in general position with respect to
the arrangement $\check A=\A\unite\tilde\A$ and use notation as in
\prettyref{sec:hyperplane}.
By \prettyref{prop:hk-barhk} and induction on the dimension of $V$
\begin{equation*}
i_!(h_k(c)\iprod \tilde h_l(d))
  =i_!(h_k(c))\iprod i_!(\tilde h_l(d))
  =h^H_{k-1}(\eta_\ast(c)) \iprod \tilde h^H_{l-1}(\eta_\ast(d))
  = 0,
\end{equation*}
since the arrangement $A^H$ is again in general position with respect
to $\tilde A^H$ and $(k-1)+(l-1)<n-1$. 
This implies 
$h_k(c)\iprod \tilde h_l(d)\in\ker i_!=\im\check h_0$
by \prettyref{prop:hk-barhk}.  But 
$h_k(c)\iprod \tilde h_l(d)\in\dsum_{j>0}\im\check h_j$ by
\prettyref{prop:h0-vanish} and hence 
$h_k(c)\iprod \tilde h_l(d)=0$.
\end{proof}

\begin{proof}[Proof of \prettyref{prop:vanishing}]\label{proof:vanishing}
We choose a neighbourhood $U$ of~$\Union P\A$ such that the
inclusion $(PV, \Union P\A)\to (PV, U)$ is a homotopy equivalence.
We then choose a copy $\tilde\A$ of $\A$ also contained in $U$, in general
position with respect to $\A$ and such that the diagram
\begin{equation*}\xymatrix{
    \C P^i\times\dpairp Qin
    \ar[d]^{f^i}
    \ar[r]^-{\tilde f^i}
    &
    (PV, \Union P\tilde\A)
    \ar[d]^{\text{incl.}}
    \\
    (PV, \Union P\A)
    \ar[r]^-{\text{incl.}}
    &
    (PV, U)
}\end{equation*}
commutes up to homotopy.
Because of the commutativity of
\begin{equation*}
\xymatrix{
    H_\ast(PV, \Union P\A)\tens H_\ast(PV, \Union P\A)
    \ar[r]^-\iprod
    &
    H_\ast(PV, \Union P\A)
    \ar[d]_{\text{incl}_\ast}^\isom
    \\
    H_\ast(\dpair Qkn)\tens H_\ast(\dpair Qln)
    \ar[u]_{h_k\tens h_l}
    \ar[d]^{h_k\tens\tilde h_l}
    &
    H_\ast(PV, U)
    \\
    H_\ast(PV, \Union P\A)\tens H_\ast(PV, \Union P\tilde\A)
    \ar[r]^-\iprod
    &
    H_\ast(PV, \Union P\A\unite\Union P\tilde\A)
    \ar[u]^{\text{incl}_\ast}
}\end{equation*}
the result follows from \prettyref{prop:vanishing-gp}
\end{proof}

\section{The main theorem}
\label{sec:main-proof}
\begin{thm}   \label{thm:main}
  Let $c\in H_\ast(Q_{[k,n]},Q_{[k,n)})$, 
    $d\in H_\ast(Q_{[l,n]},Q_{[l,n)})$.
Then
\begin{equation*}
 h_k(c)\iprod h_l(d)=
 h_{k+l-n}(c\mintimes d)
\end{equation*}
with the convention that $h_i=0$ for $i<0$.
\end{thm}

\begin{proof}
We proceed by induction on the dimension of~$V$.

For $k+l<n$ the conjecture is covered by \prettyref{prop:vanishing}.
For $k+l\ge n$ and $\dim V=1$ it is covered by 
\prettyref{prop:graded-product} (and trivial anyway).
For $k+l\ge n$ and $\dim V>1$ we choose a hyperplane in general
position with respect to the arrangement and adopt the notation of
\prettyref{sec:hyperplane}.  By induction and
\prettyref{prop:hk-barhk} (for $k+l>n$) or \prettyref{prop:vanishing} 
(for $k+l=n$)
\begin{multline*}
i_!(h_k(c)\iprod h_l(d))
  =i_!(h_k(c))\iprod i_!(h_l(d))
  =h^H_{k-1}(\eta_\ast(c)) \iprod h^H_{l-1}(\eta_\ast(d))
  =\\=h^H_{k+l-n-1}(\eta_\ast(c)\mintimes \eta_\ast(d))
  =h^H_{k+l-n-1}(\eta_\ast(c\mintimes d))
  =i_!(h_{k+l-n}(c\mintimes d)).
\end{multline*}
Again by \prettyref{prop:hk-barhk}, this implies 
$h_k(c)\iprod h_l(d)-h_{k+l-n}(c\mintimes d)\in\ker i_!=\im h_0$.  
But $h_k(c)\iprod h_l(d)-h_{k+l-n}(c\mintimes d)\in
\dsum_{i>k+l-n}\im h_i$ by \prettyref{prop:graded-product}. Therefore 
$h_k(c)\iprod h_l(d)-h_{k+l-n}(c\mintimes d)=0$.
\end{proof}

\begin{rem} \label{rem:affine-products}
Continuing \prettyref{rem:affine-homology}, this yields a
description of the cohomology ring of an affine complex arrangement
with intersection poset~$Q'$.  All of the isomorphisms in that remark
are induced by inclusions and therefore respect the
product~$\mintimes$.  Denoting, for $u\in Q'_{[0,n]}$, the map
$H(\Delta[u,V],\Delta[u,V)\unite\Delta(u,V])\to H(PV,\Union P\A)$
arising from these isomorphisms and $h_{d(u)}$ by $h_u$ we
obtain
\begin{equation*}
  h_u(c)\iprod h_v(d)=
  \begin{cases}
     h_{u\lmin v}(c\mintimes d),&
        \text{if $u\lmin v\in Q'_{[0,n]}$}
        \\&\quad\text{and $d(u\lmin v)=d(u)+d(v)-n$,}\\
     0,&\text{otherwise,}
  \end{cases}
\end{equation*}
for $c\in H_\ast(\Delta[u,V],\Delta[u,V)\unite\Delta(u,V])$, 
$d\in H_\ast(\Delta[v,V],\Delta[v,V)\unite\Delta(v,V])$.
\end{rem}

\section{Projective \texorpdfstring{$c$}{c}-arrangements}
\let\seci\subsection
\label{sec:c-arr}
For special classes of arrangements, one can hope to derive a simpler
description of the cohomology ring of the complement from
\prettyref{thm:main}.  The case that is probably easiest to handle is
that of $c$-arrangements.
\begin{defn}\label{def:c-arr}
For a positive integer $c$, we call $\A$ a $c$-arrangement, if every
$A\in\A$ is a subspace of codimension~$c$ and $d(q)$ is an integral
multiple of~$c$ for every~$q\in Q$.
\end{defn}

\seci{A presentation of the cohomology ring}

\begin{defn}\label{def:dependent}
We call a subset $M$ of~$\A$ independent, if 
$n-d\left(\Intersection M\right) =\sum_{A\in M}(n- d(A))$, 
dependent, if it is not independent, and minimally dependent, if it
is dependent but all of its proper subsets are independent.
\end{defn}

In this section we will prove:

\begin{thm}\label{thm:c-arr}
Let $\A$ be a $c$-arrangement, $\card\A-1\eqd t\ge0$,
$\A=\set{A_0,\dots,A_t}$.  Let $R$ be the free graded commutative (in
the graded sense) ring over the set of generators $\set
x\unite\set{y_i\colon 1\le i\le t}$ with $\adeg x=2$,
$\adeg{y_i}=2c-1$.  Let $I$ be the ideal generated by
\begin{align*}
&\set{\sum_{j=0}^r (-1)^j y_{i_0}\cdots\hat y_{i_j}\cdots y_{i_r} 
       \colon \text{$i_0<\cdots<i_r$, $\set{A_{i_j}}$ is minimally dependent.}}
  \\&\unite
  \set{y_{i_1}\cdots y_{i_r}\colon
       \text{$i_1<\cdots<i_r$,
	     $\set{A_0}\unite\set{A_{i_j}}$ is minimally dependent.
            }
      }\\&\unite\set{x^c}.
\end{align*}
The map
\begin{align*}
\pi\colon
    R&\to H^\ast\left(PV\wo \Union P\A\right), \\
    x&\mapsto P(h_{n-1}([\simp{V}])),\\
    y_i&\mapsto P(h_{n-c}([\simp{A_i,V}-\simp{A_0,V}])),
\end{align*}
where 
$P\colon H_\ast(PV, \Union P\A)\xto\isom 
H^\ast\left(PV\wo \Union P\A\right)$
denotes Poincar\'e duality,
is an epimorphism and $\ker\pi=I$.
\end{thm}

We now fix the arrangement $\A=\set{A_0,\dots,A_t}$.

\begin{rem}
For $c=1$ the complement $PV\wo\Union P\A$ can, as in
\prettyref{rem:affine-homology}, be regarded as the
complement in the affine space $PV\wo PA_0$ of the linear hyperplane
arrangement $\A'\deq\set{PA_i\wo PA_0\colon 1\le i\le t}$.  In this case, the
generator $x$ and the corresponding relation can be omitted.

If $A_0$ is in general position with respect to $\A\wo\set{A_0}$, the
second kind of generators does not occur.  This is in particular the
case if the arrangement $\A'$ is central, i.e.\ if
$\Intersection\A'\ne\emptyset$. In this case the theorem reduces to
the description of the cohomology ring of the complement of~$\A'$ 
given by Orlik and Solomon\cite{orliksolomon80}.
\end{rem}
\seci{The atomic complex}
We now turn to the proof of the theorem.  When using simplicial chain
complexes, we will always use the complex of non-degenerate simplices
and view it as the complex of all simplices modulo degenerate
simplices if necessary.

\begin{defn}
For an integer~$k$ with $0\le k\le n$, we define $S_k$ to be the
simplicial complex which has the vertex set $\set{0,\dots,t}$ and as simplices
the sets $I\subset\set{0,\dots,t}$ with
$d\left(\Intersection_{i\in I} A_i\right)\ge k$.  This is the \emph{atomic
complex} of $Q_{[k,n]}$.  We also define $D^k$ to
be the reduced ordered (using the natural order of $\set{0,\dots,t}$) 
simplicial chain complex of $S_k$ shifted by one, 
i.e.\ $D^k_r=\tilde C_{r-1}(S_k)$ and
in particular $D^k_0\isom\Z$ generated by the empty simplex.
\end{defn}
As is well known, the atomic complex and the order complex, of
$Q_{[k,n)}$ in this case, are homotopy equivalent.  We describe a
homotopy equivalence to fix a concrete isomorphism between their
homology groups.  Before doing this, we state a useful lemma.
\begin{lem}
Let $P_0$, $P_1$ be posets, $P'_i\subset P_i$.  If $f,g\colon
(P_0,P_0')\to (P_1,P_1')$ are order preserving functions such that
$f(p)\le g(p)$ for all $p\in P_0$, then the maps $f,g\colon(\Delta
P_0,\Delta P_0')\to(\Delta P_1,\Delta P_1')$ are homotopic.
\end{lem}

\begin{proof}
The map $H\colon\set{0,1}\times P_0\to P_1$ defined by $H(0,x)\deq
f(x)$, $H(1,x)\deq g(x)$ is order preserving and hence yields the desired
homotopy
\begin{equation*}
I\times (\Delta P_0,\Delta P_0')\homeo
\left(\Delta (\set{0,1}\times P_0),\Delta (\set{0,1}\times P_0')\right)
\xto{\,H\,}(\Delta P_1,\Delta P_1'),
\end{equation*}
 where
we view $\set{0,1}$ as a poset.
\end{proof}

\begin{rem}
This lemma is a special case of
\cite[Prop.~2.1]{segal-classifying-spaces} which is proved in the same way.
\end{rem}

\begin{defprop}\label{prop:s-homot-equiv}
We denote the face poset of $S_k$ by $FS_k$, but order it by $M\le M'$
if $M'$ is a face of $M$, that is if $M'\subset M$.  We also set
$\tilde FS_k\deq FS_k\unite\set\emptyset$.
The map
\begin{align*}
  s\colon (\tilde FS_k,FS_k)&\to (Q_{[k,n]},Q_{[k,n)})\\
    M&\mapsto\Intersection \set{A_i\colon i\in M}
\end{align*}
is then order preserving and moreover satisfies $s(M\lmin
M')=s(M\unite M')=s(M)\intersect s(M')=s(M)\lmin s(M')$, if one side,
and therefore the other, exists.

With these definitions, the map $s\colon(\Delta \tilde FS_k,\Delta
FS_k)\to\dpairp Qkn$ is a homotopy equivalence.
\end{defprop}

\begin{rem}
$\Delta FS_k$ is the barycentric subdivision of $S_k$, 
and $\Delta\tilde FS_k$ is a cone over $\Delta FS_k$.
\end{rem}

\begin{proof}
We define an order preserving map
\begin{align*}
  r\colon (Q_{[k,n]},Q_{[k,n)})&\to(\tilde FS_k,FS_k),\\
    q&\mapsto\set{i\colon A_i\supset q}.
\end{align*}
We have $s(r(q))\ge q$ for $q\in Q_{[k,n]}$ and $r(s(i))\le i$ for
$i\in\tilde FS_k$.  Hence, by the preceding lemma $r$~is a homotopy
inverse to $s$, when both maps are regarded as simplicial maps between
order complexes.
\end{proof}

\begin{defprop}
We define chain maps
\begin{align*}
f^k\colon D^k_r&\to C_r(\dpair Qkn)\\
  \simp{i_1,\dots,i_r}&\mapsto 
          \simp{A_{i_1},V}\mintimes\cdots\mintimes\simp{A_{i_r},V}.
\end{align*}
For $r=0$ this is to be understood as 
$f^k(\simp{})=\simp{V}$.
\end{defprop}

\begin{nota}
To simplify the following calculations, we set set
\begin{equation*}
\alpha_i\deq\simp{A_i,V}\in C_1(\Delta Q_{[n-c,n]})
\end{equation*}
and sometimes write the multiplication $\mintimes$ as juxtaposition.
\end{nota}

\begin{proof}
To see that $f^k$ is well-defined, we have to check that the right
hand side is in $C_\ast(\Delta Q_{[k,n]})$.  But 
$d\left(A_{i_1}\lmin\cdots\lmin A_{i_r}\right)\ge k$ by definition of
$S_k$ and hence $D^k_r$.

To see that $f^k$ is a chain map, we calculate
\begin{equation*}
\begin{split}
\dd (f^k(\simp{i_1,\dots,i_r}))
    &= \sum_{j=1}^r (-1)^{j+1}
       \simp{A_{i_1},V}\mintimes\cdots
       \mintimes\dd\simp{A_{i_j},V}\mintimes
       \cdots\mintimes\simp{A_{i_r},V}
  \\&= \sum_{j=1}^r (-1)^{j+1}\alpha_{i_1}\cdots\alpha_{i_{j-1}}
                              \simp V
			      \alpha_{i_{j+1}}\cdots\alpha_{i_r}
    \\&\qquad-\sum_{j=1}^r (-1)^{j+1}\alpha_{i_1}\cdots\alpha_{i_{j-1}}
                              \simp{A_{i_j}}
			      \alpha_{i_{j+1}}\cdots\alpha_{i_r}
  \\&= \sum_{j=1}^r (-1)^{j+1}\alpha_{i_1}\cdots\alpha_{i_{j-1}}
                              \hat\alpha_{i_j}
			      \alpha_{i_{j+1}}\cdots\alpha_{i_r}
    \\&\qquad-\sum_{j=1}^r (-1)^{j+1}\alpha_{i_1}\cdots\alpha_{i_{j-1}}
                              \simp{A_{i_j}}
			      \alpha_{i_{j+1}}\cdots\alpha_{i_r}.
\end{split}
\end{equation*}
The first summand equals $f^k(\dd\simp{i_1,\dots,i_r})$. 
Since 
$\simp{A_i}\in\C_0(Q_{[n-c,n)})$,
the second
summand is in 
$C_\ast(\Delta Q_{[k,n-c]})\subset C_\ast(\Delta Q_{[k,n)})$.
\end{proof}

\begin{prop} \label{prop:fk-iso}
The induced maps $f^k_\ast\colon H(D^k)\to H(\dpair Qkn)$ are isomorphisms.
\end{prop}

\begin{proof}
Defining
\begin{align*}
\bar f^k\colon D^k_r&\to C_r(\Delta \tilde FS_k,\Delta FS_k)\\
  \simp{i_1,\dots,i_r}&\mapsto 
          \simp{\set{i_1},\emptyset}\mintimes
          \cdots\mintimes\simp{\set{i_r},\emptyset}
\end{align*}
the diagram
\begin{equation*}\xymatrix{
    &
    H_r(D^k)
    \ar[d]_-{\bar f^k_\ast}
    \ar[ld]_{f^k_\ast}
    \ar[rd]^{\mathop{sd}_\ast}_\isom
    \\
    H_r(\dpair Qkn)
    &
    H_r(\Delta \tilde FS_k,\Delta FS_k)
    \ar[l]_-{s_\ast}^-\isom
    \ar[r]^-\dd_-\isom
    &
    \tilde H_{r-1}(\Delta FS_k)
}\end{equation*}
commutes, where $\mathop{sd}$ is the barycentric subdivision map
$\tilde C_\ast(S_k)\to\tilde C_\ast(\Delta FS_k)$.
The connecting homomorphism is an
isomorphism, because $\tilde FS_k$ has the maximum~$\emptyset$. 
The map~$s_\ast$
is an isomorphism because of \prettyref{prop:s-homot-equiv}.  
It follows that $f^k_\ast$ is an isomorphism.
\end{proof}
\seci{Proof of the presentation}
The chain maps $f^k$ would be more useful in a situation in which the
chains $\simp{A_i,V}$ are cycles. For example, think of affine
arrangements, where, acoording to \prettyref{rem:affine-homology},
 $(\Delta Q_{[k,n]},\Delta Q_{[k,n)}\unite\Delta
Q_{(k,n]})$ takes the place of $\dpairp Qkn$.  In our situation they
are not.  The chains $\simp{A_i,V}-\simp{A_j,V}$ however are cycles,
we will therefore replace the maps $f^k$ by the following maps.
\begin{defprop}
For a $c$-arrangement $\A$, we define chain maps
\begin{align*}
g^k\colon D^k_r&\to C_r(\dpair Qkn)\\
  \simp{i_1,\dots,i_r}&\mapsto 
  \begin{cases}
    (\simp{A_{i_1},V}-\simp{A_0,V})
      \mintimes\cdots\mintimes
    (\simp{A_{i_r},V}-\simp{A_0,V}),
      & r=a,\\
    0,&r\ne a,
  \end{cases}
\end{align*}
where $a$ is defined by $n-(a+1)c<k\le n-ac$.
\end{defprop}

\begin{proof}We check that $g^k$ is a well-defined chain map.
For $r>a$ we have $C_r(\dpair Qkn)\isom0$, since 
$n-k<(a+1)c\le rc$.  So we just have to show that 
$g^k(\simp{i_1,\dots,i_a})$ is a cycle in $C_a(\dpair Qkn)$.
This is true, because each $\simp{A_i,V}-\simp{A_0,V}$ is a cycle in 
$C_1(\dpair Q{n-c}n)$ and $n-ac\ge k$.
\end{proof}

\begin{prop} \label{prop:fk-gk-homot}
The maps $f^k$ and $g^k$ are chain homotopic.
\end{prop}

\begin{proof}
We define
\begin{align*}
  K\colon D^k_r&\to C_{r+1}(\dpair Qkn)\\
\simp{i_1,\dots,i_r}&\mapsto 
  \begin{cases}
    f^k\left(\simp{0,i_1,\dots,i_r}\right),
      & r<a,\\
    0,&r\ge a.
  \end{cases}
\end{align*}
The right hand side is well defined, because for $r<a$ we have
\[d(A_0\intersect A_{i_1}\intersect\dots\intersect A_{i_r})
  \ge n-(r+1)c\ge n-ac\ge k.\]
We calculate $K\dd+\dd K$.

For $r<a$:
\begin{align*}
  &(K\dd+\dd K)\simp{i_1,\dots,i_r}
  \\&\quad=f^k\left(\sum_{j=1}^r (-1)^{j+1}
                                 \simp{0,i_1,\dots,\hat i_j,\dots,i_r}
               \right)
      +\dd f^k\left(\simp{0,i_1,\dots,i_r}\right)
  \\&\quad=f^k\left(\sum_{j=1}^r (-1)^{j+1}
                                 \simp{0,i_1,\dots,\hat i_j,\dots,i_r}
                    + \dd\simp{0,i_1,\dots,i_r}
              \right)
  \\&\quad=f^k(\simp{i_1,\dots,i_r})
          =(f^k-g^k)\simp{i_1,\dots,i_r}.
\end{align*}

For $r=a$:
We first calculate
\begin{align*}
g^k(\simp{i_1,\dots,i_a})
  &=(\alpha_{i_i}-\alpha_0)\cdots(\alpha_{i_a}-\alpha_0)
  \\&=\alpha_{i_1}\cdots\alpha_{i_a}
      -\sum_{j=1}^a \alpha_{i_1}\cdots\alpha_{i_{j-1}}
                    \alpha_0
		    \alpha_{i_{j+1}}\cdots\alpha_{i_a}
  \\&=\alpha_{i_1}\cdots\alpha_{i_a}
      +\sum_{j=1}^a (-1)^j \alpha_0
                    \alpha_{i_1}\cdots
                    \hat\alpha_{i_j}
		    \cdots\alpha_{i_a}
\end{align*}
and with this
\begin{align*}
(K\dd+\dd K)\simp{i_1,\dots,i_a}
  &= f^k\left(\sum_{j=1}^a (-1)^{j+1}\simp{0,i_1,\dots,\hat i_j,\dots,i_a}
        \right)
  \\&= \sum_{j=1}^a (-1)^{j+1} \alpha_0
                    \alpha_{i_1}\cdots
                    \hat\alpha_{i_j}
		    \cdots\alpha_{i_a}
  \\&= (f^k-g^k)\simp{i_1,\dots,i_a}.
\end{align*}

For $r>a$ we have 
$(K\dd+\dd K)\simp{i_1,\dots,i_r}=0=(f^k-g^k)\simp{i_1,\dots,i_r}$,
since
$C_r(\dpair Qkn)\isom0$ as noted before.
\end{proof}

\begin{prop}
The map $\pi$ is surjective.
\end{prop}

\begin{proof}
By \prettyref{prop:fk-iso}, \prettyref{prop:fk-gk-homot}, and of
course \prettyref{prop:proj-homology}, 
$H^\ast(PV\wo\Union P\A)$ is additively generated by the elements 
$P\left(h_k\left([g^k\simp{i_i,\dots,i_r}]\right)\right)$ with
$k\le n-rc$.  By \prettyref{thm:main} 
\begin{align*}
&P\left(h_k\left([g^k\simp{i_i,\dots,i_r}]\right)\right)
  \\&\quad=
  P(h_k([(\alpha_{i_1}-\alpha_0)\mintimes\cdots\mintimes
	 (\alpha_{i_r}-\alpha_0)]))
  \\&\quad=
  P(h_{n-1}([\simp V]))^{n-k-rc}
      P(h_{n-c}([\alpha_{i_1}-\alpha_0]))\cdots
      P(h_{n-c}([\alpha_{i_r}-\alpha_0]))
  \\&\quad=
  \pi(x)^{n-k-rc}\pi(y_{i_1})\cdots\pi(y_{i_r}).
\end{align*}
This shows that $\pi$ is surjective.  
\end{proof}

\begin{prop}
$I\subset\ker\pi$.
\end{prop}

\begin{proof}
First of all
\begin{align*}
\pi(x^c) 
  &=
  P(h_{n-1}([\simp V]))^c
  \\&=P(h_{n-c}([\simp V]))
  \\&=P(h_{n-c}([\dd\simp{A_0,V}]))=P(h_{n-c}(0))=0.
\end{align*}
If $\set{A_{i_0},\dots,A_{i_r}}$ is minimally dependent, then 
$d\left(\Intersection_j A_{i_j}\right)=n-rc$ and
\begin{align*}
  0&=P(h_{n-rc}(g^{n-rc}_\ast([\dd\simp{i_0,\dots,i_r}])))
  \\&=(P\cmps h_{n-rc}\cmps g^{n-rc}_\ast)
      \left(\left[
        \sum_{j=0}^r (-1)^j\simp{i_0,\dots,\hat i_{j},\dots,i_r}
      \right]\right)
  \\&=\pi\left(\sum_{j=0}^r (-1)^j y_{i_0}
				  \cdots\hat y_{i_j}
				  \cdots y_{i_r}
        \right)
\end{align*}
and similarly if $\set{A_0,A_{i_1},\dots,A_{i_r}}$ is minimally
dependent, then
\begin{align*}
  0&=P(h_{n-rc}(g^{n-rc}_\ast([\dd\simp{0,i_1,\dots,i_r}])))
  \\&=(P\cmps h_{n-rc}\cmps g^{n-rc}_\ast)
      \left(\left[
	\simp{i_1,\dots,i_r}+
        \sum_{j=1}^r (-1)^j\simp{0,i_1,\dots,\hat i_{j},\dots,i_r}
      \right]\right)
  \\&=\pi\left(y_{i_1}\cdots y_{i_r}
         \right)
\end{align*}
as claimed.
\end{proof}

\begin{lem} \label{lem:dep-in-I}
If $i_0<\cdots<i_r$ and $\set{A_{i_j}}$ is dependent, then
$y_{i_0}\cdots y_{i_r}\in I$ and 
$\sum_j (-1)^j y_{i_0}\cdots \hat y_{i_j}\cdots y_{i_r}\in I$.
\end{lem}

\begin{proof}
Let $\set{A_{i_0},\dots,A_{i_r}}$ be dependent.  To show 
$y_{i_0}\cdots y_{i_r}\in I$ we may assume that the set is
minimally dependent.  Then
$\sum_j (-1)^j y_{i_0}\cdots\hat y_{i_j}\cdots y_{i_r}\in I$ and
$y_{i_0}
 \left(\sum_j (-1)^j y_{i_0}\cdots\hat y_{i_j}\cdots y_{i_r}
 \right)=y_{i_0}\cdots y_{i_r}$ since $y_0^2=0$.

For the second part of the lemma we may assume that 
$\set{A_{i_j}\colon j\le s}$ is minimally
dependent.  Then 
\begin{multline*}
\sum_j (-1)^j y_{i_0}\cdots \hat y_{i_j}\cdots y_{i_r}
  =
  \\\underbrace{\left(\sum_{j=0}^s (-1)^j 
              y_{i_0}\cdots\hat y_{i_j}\cdots y_{i_s}
   \right)}_{\in I}
   y_{i_{s+1}}\cdots y_{i_r}
  \\+\underbrace{y_{i_0}\cdots y_{i_s}}_{\in I}
   \sum_{j=s+1}^r (-1)^j y_{i_{s+1}}\cdots\hat y_{i_j}\cdots y_{i_r}
  \in I
\end{multline*}
as claimed.
\end{proof}

\begin{prop}
$\ker\pi\subset I$.
\end{prop}

\begin{proof}
Let $z\in\ker\pi$.  We want to show $z\in I$.  We may assume that z is
a linear combination of elements $x^sy_{i_1}\cdots y_{i_r}$ with 
$0\le s<c$, $i_1<\cdots<i_r$ and $\set{A_{i_1},\dots,A_{i_r}}$
independent.  Since
$\pi(x^sy_{i_1}\cdots y_{i_r})\in\im(P\cmps h_{n-cr-s})$ and $r$
and~$s$ are determined by $cr+s$, we may assume $z$ to be homogenous
in $r$ and $s$, i.e.\ 
$z=x^s\sum_{i_1<\cdots<i_r}\lambda_i y_{i_1}\cdots y_{i_r}$.
We set $k\deq n-cr-s$.
The chain 
\begin{equation*}
z'\deq\sum_{i_1<\cdots<i_r}\lambda_i\simp{i_1,\dots,i_r}
      +\sum_{i_1<\cdots<i_r}\lambda_i\sum_{j=1}^r(-1)^j
                            \simp{0,i_1,\dots,\hat i_j,\dots,i_r}
\end{equation*}
is a cycle in $D^k_r$ (the second summand is a cone over the boundary
of the first summand), 
$0=\pi(z)=(P\cmps h_k)\left(g^k_\ast\left([z']\right)\right)$, and
therefore $[z']=0$ by \prettyref{prop:fk-iso} and 
\prettyref{prop:fk-gk-homot}, i.e.\ $z'$ is a boundary in $D^k$, 
which means that there exist 
$\mu_i$, $\nu_i$ such that
\[z'=\dd\left(
    \sum_{i_1<\cdots<i_r}\mu_i\simp{0,i_1,\dots,i_r}
    +\sum_{i_0<\cdots<i_r}\nu_i\simp{i_0,\dots,i_r}
  \right)
\]
and with $d\left(\Intersection_j A_{i_j}\intersect A_0\right)\ge k$
and therefore $\set{A_0}\unite\set{A_{i_j}}$ dependent for $\mu_i\ne0$
and $\set{A_{i_j}}$ dependent for $\nu_i\ne0$.  Comparing coefficients
and sorting the simplices by whether the first vertex is $0$ yields
\[
z=\sum_i\mu_i y_{i_1}\cdots y_{i_r}
 +\sum_i\nu_i\sum_j (-1)^j y_{i_0}\cdots\hat y_{i_j}\cdots y_{i_r}
 \in I
\]
as claimed.
\end{proof}

This completes the proof of \prettyref{thm:c-arr}\qed

\bibliographystyle{amsalpha}
\bibliography{math,topology,combi}

\end{document}